\documentclass[12pt]{article}
\usepackage{amssymb,amsmath, amsthm, mathtools}
\usepackage{colordvi}
\usepackage{graphicx}
\usepackage{verbatim}
\usepackage{cancel}
\usepackage{caption}

\usepackage[includeheadfoot,left=2.5cm,right=2.5cm,top=2.5cm,bottom=2.5cm,head=15pt,
bindingoffset=0pt]{geometry}
\usepackage[all,  colour, line, dvips]{xy}
{\numberwithin{equation}{section}}

\newenvironment{prf}[1][Proof]{\begin{proof}[{\bf#1. }]} {\end{proof}}  

\theoremstyle{plain}
\newtheorem{thm}[equation]{Theorem}
\newtheorem{lem}[equation]{Lemma}
\newtheorem{cor}[equation]{Corollary}

\newtheorem{prop}[equation]{Proposition}
\newtheorem*{thm*}{Theorem}

\theoremstyle{definition}

\newtheorem{defin}[equation]{Definition}
\newtheorem{rmk}[equation]{Remark}
\newtheorem{ex}[equation]{Example}

\theoremstyle{plain}
\newcounter{introthm}[section]
\newtheorem{manualnothmcanva}[introthm]{Theorem}
\newenvironment{manualnothm}[1]{\renewcommand{\theintrothm}{#1}
\begin{manualnothmcanva}} {\end{manualnothmcanva}}

\newtheorem*{ex*}{Example}

\def\myVSPACEfigure{\vspace{0.7cm}}
\def\git{/\!\!/}
\def\biggit{\Bigr/\!\!\!\Bigl/}
\def\cal{\mathcal}

\def\ccS{\mathcal{S}}

\def\conv{\operatorname{conv}}

\def\SS{\mathbb{S}}
\def\boldS{\mathbb{S}}

\def\half{\frac{1}{2}}

\def\CC{{\mathbb C}}

\def\NN{{\mathbb N}}

\def\PP{{\mathbb P}}

\def\SS{{\mathbb S}}
\def\RR{{\mathbb R}}
\def\TT{{\mathbb T}}
\def\ZZ{{\mathbb Z}}

\def\ra{\rightarrow}

\def\mono{\hookrightarrow}
\def\epi{\twoheadrightarrow}

\def\ccB{{\cal B}}

\def\ccE{{\cal E}}

\def\ccG{{\cal G}}

\def\ccL{{\cal L}}

\def\ccO{{\cal O}}
\def\ccP{{\cal P}}

\def\ccS{{\cal S}}
\def\ccT{{\cal T}}

\def\ccV{{\cal V}}
\def\ccX{{\cal X}}
\def\ccY{{\cal Y}}

\def\hM{{\widehat M}}

\def\Mg{M^{gr}}

\def\Proj{\operatorname{Proj}}
\def\Spec{\operatorname{Spec}}

\def\Hom{\operatorname{Hom}}

\def\deg{\operatorname{deg}}
\def\id{\operatorname{id}}

\def\comp{\operatorname{comp}}

\newcommand\rootthree{1.73205}
\newcommand\fiverootthree{8.66025403784439}

\def\CHERRY{
\POS( 0, 0) ; \POS( 4, 6 ) **\crv{~**\dir{-} ( 4, 2 )};
\POS( 4, 6) ; \POS( 0, 10) **\crv{~**\dir{-} ( 4, 10)};
\POS(-4, 6) ; \POS( 0, 10) **\crv{~**\dir{-} (-4, 10)};
\POS( 0, 0) ; \POS(-4, 6 ) **\crv{~**\dir{-} (-4, 2 )}; 
}

\def\BALLOON#1#2{
\begin{xy}
(0,0);(0.3,0):;
\POS( 0, 0) ; \POS( 4, 6 ) **\crv{~**\dir{#1} ( 4, 2 )};
\POS( 4, 6) ; \POS( 0, 10) **\crv{~**\dir{#1} ( 4, 10)};
\POS(-4, 6) ; \POS( 0, 10) **\crv{~**\dir{#1} (-4, 10)};
\POS( 0, 0) ; \POS(-4, 6 ) **\crv{~**\dir{#1} (-4, 2 )}; 
\POS( 0, 0) ; \POS( 0, -6) **\dir{#2};
\end{xy}
}

\def\BALLOONfont{
\begin{xy}\save
< 0mm,0.2mm>;<0.1666mm,0.2mm>:; 
\drop{\CHERRY}
\POS( 0, 0.5) ; \POS( 0,-8 ) **@{-};
(-5,0)*{}; (5,0)*{};
\restore\end{xy}}

\def\BALLOONscriptfont{
\begin{xy}\save
\POS<0mm,0mm>;\POS<0.1mm,0mm>:;
\drop{\CHERRY}
\POS( 0, 0.2) ; \POS( 0,-10 ) **@{-};
(-5,0)*{}; (5,0)*{};
\restore\end{xy}}

\def\TRIPOD {    
\begin{xy}
=(              0,  0) "0";
=(              0, 10) "A"; 
=(-\fiverootthree, -5) "B"; 
=( \fiverootthree, -5) "C";
<0mm,-0.1mm>;<0.2mm,-0.1mm>:; 
"0" ; "A" **@{-};
"0" ; "B" **@{-};
"0" ; "C" **@{-};
\end{xy}}

\def\LMfont{
\begin{xy}\save
< 0mm,0.1666mm>;<0.15mm,0.1666mm>:; 
\drop{\CHERRY}
\POS( 0, 0.6) ; \POS( 0,-6 ) **@{-};
\POS( 0,-6) ; \POS(-9,-12) **@{-};
\POS( 0,-6) ; \POS( 9,-12) **@{-};
\restore\end{xy}}

\def\LMscriptfont{
\begin{xy}\save
< 0mm,0mm>;<0.1mm,0mm>:; 
\drop{\CHERRY}
\POS( 0, 0.6) ; \POS( 0,-6 ) **@{-};
\POS( 0,-6) ; \POS(-9,-12) **@{-};
\POS( 0,-6) ; \POS( 9,-12) **@{-};
\restore\end{xy}}

\def\HAMMOCKfont{
\begin{xy}\save
(-2,0)*{} ; (-1,0)*{} **@{-};
(-1,0)*{} ; ( 1,0)*{}  **\crv{~**\dir{-} (0,1.75)};
(-1,0)*{} ; ( 1,0)*{}  **\crv{~**\dir{-} (0,-1.75)};
( 2,0)*{} ; ( 1,0)*{} **@{-};
\restore\end{xy}}

\def\DUMBBELLfont{
\begin{xy}
\POS<0.49mm,0.4mm>;\POS<0.6mm,0.4mm>:;
\drop{\CHERRY}
\drop{\POS( 0, 0); \POS( 0,-6 )**@{-};}
\POS( 0,-3); \POS(-1,-3):;
\drop{\CHERRY}
\POS(-8,0)*{};
\end{xy}}

\def\DUMBBELLscriptfont{
\begin{xy}
\POS<0.36mm,0.09mm>;\POS<0.43mm,0.09mm>:;
\drop{\CHERRY}
\drop{\POS( 0, 0); \POS( 0,-6 )**@{-};}
\POS( 0,-3); \POS(-1,-3):;
\drop{\CHERRY}
\POS(-8,0)*{};
\end{xy}}

\def\THETAfont{
\begin{xy}
< 0mm,0.37mm>;<0.76mm,0.37mm>:;
(0,0);(3  , 0) **@{-};
(0,0);(1.5, 2) **\crv{~**\dir{-}(0,2)} ;
(3,0);(1.5, 2) **\crv{~**\dir{-}(3,2)} ;
(0,0);(1.5,-2) **\crv{~**\dir{-}(0,-2)} ;
(3,0);(1.5,-2) **\crv{~**\dir{-}(3,-2)} ;
\end{xy}
}

\def\THETAscriptfont{
\begin{xy}
< 0mm,0mm>;<0.5mm,0mm>:;
(0,0);(3  , 0) **@{-};
(0,0);(1.5, 2) **\crv{~**\dir{-}(0,2)} ;
(3,0);(1.5, 2) **\crv{~**\dir{-}(3,2)} ;
(0,0);(1.5,-2) **\crv{~**\dir{-}(0,-2)} ;
(3,0);(1.5,-2) **\crv{~**\dir{-}(3,-2)} ;
\end{xy}
}
\def\TREEfont{
\begin{xy}\save
< 0mm,0.1666mm>;<0.1666mm,0.1666mm>:;
(0, 3) ; (-9, 9 ) **@{-};
(0, 3) ; ( 9, 9 ) **@{-};
(0, 3) ; ( 0,-6 ) **@{-};
(0,-6) ; (-9,-12) **@{-};
(0,-6) ; ( 9,-12) **@{-};
\restore\end{xy}}

\def\TREEbackup#1#2#3#4#5[#6][#7][#8][#9]{ 
\begin{xy} 
\save
(0, 0)*{} ; (-2, 2)*{} **@{#1}, (-3, 2)  *={\scriptstyle{#6}};
(0, 0)*{} ; ( 2, 2)*{} **@{#2}, ( 3, 2)  *={\scriptstyle{#7}};
(0, 0)*{} ; ( 0,-3)*{} **@{#3}; 
(0,-3)*{} ; (-2,-5)*{} **@{#4}, (-3,-5)  *={\scriptstyle{#8}};
(0,-3)*{} ; ( 2,-5)*{} **@{#5}, ( 3,-5)  *={\scriptstyle{#9}};
\restore
\end{xy}
}

\def\TREE#1#2#3#4#5[#6][#7][#8][#9]{ 
\begin{xy} 
\save
(0, 1.5);(1.43,1.5):;
(0, 0)*{} ; (-2, 2)*{} **@{#1}, (-3, 2)  *={\scriptstyle{#6}};
(0, 0)*{} ; ( 2, 2)*{} **@{#2}, ( 3, 2)  *={\scriptstyle{#7}};
(0, 0)*{} ; ( 0,-3)*{} **@{#3}; 
(0,-3)*{} ; (-2,-5)*{} **@{#4}, (-3,-5)  *={\scriptstyle{#8}};
(0,-3)*{} ; ( 2,-5)*{} **@{#5}, ( 3,-5)  *={\scriptstyle{#9}};
\restore
\end{xy}
}

\def\LM#1#2#3#4[#5][#6]{ 
\begin{xy}
\save
( 0, 0)*{} ; ( 2, 3)*{} **\crv{~**\dir{#1} (2,1)};
( 2, 3)*{} ; ( 0, 5)*{} **\crv{~**\dir{#1} (2,5)};
(-2, 3)*{} ; ( 0, 5)*{} **\crv{~**\dir{#1} (-2,5)};
( 0, 0)*{} ; (-2, 3)*{} **\crv{~**\dir{#1} (-2,1)};
( 0, 0)*{} ; ( 0,-3) **@{#2};
( 0,-3)*{} ; (-2,-5)*{} **@{#3}, (-3,-5)*={\scriptstyle{#5}};
( 0,-3)*{} ; ( 2,-5)*{} **@{#4}, (3,-5)*={\scriptstyle{#6}};
\restore
\end{xy}}

\def\HAMMOCK#1#2#3#4[#5][#6] { 
\begin{xy}
\save
(0,-5)*{} ; (0,-2)*{} **@{#1},  (1, -5)*={\scriptstyle{#5}};
(0,-2)*{} ; (0, 2)*{}  **\crv{~**\dir{#2} (-3,0)};
(0,-2)*{} ; (0, 2)*{}  **\crv{~**\dir{#3} ( 3,0)};
(0, 5)*{} ; (0, 2)*{} **@{#4},  (1, 5)*={\scriptstyle{#6}};
\restore
\end{xy} }

\def\HAMMOCKbak#1#2#3#4[#5][#6] { 
\begin{xy}
\save
(-6,0)*{} ; (-3,0)*{} **@{#1},  (-7,0)*={\scriptstyle{#5}};
(-3,0)*{} ; (3,0)*{}  **\crv{~**\dir{#2} (0,5)};
(-3,0)*{} ; (3,0)*{}  **\crv{~**\dir{#3} (0,-5)};
(6,0)*{} ; (3,0)*{} **@{#4},  (7,0)*={\scriptstyle{#6}};
\restore
\end{xy} }

\def\LEGS#1#2#3#4#5#6{
\save
"A" ; "A"+"A"*{} **@{#1};
"B" ; "B"+"B"*{} **@{#2};
"C" ; "C"+"C"*{} **@{#3};
"D" ; "D"+"D"*{} **@{#4};
"E" ; "E"+"E"*{} **@{#5};
"F" ; "F"+"F"*{} **@{#6};
\restore
}

\def\HEX#1#2#3#4#5#6{
\save
"A" ; "B" **@{#1};
"B" ; "C" **@{#2};
"C" ; "D" **@{#3};
"D" ; "E" **@{#4};
"E" ; "F" **@{#5};
"F" ; "A" **@{#6};
\restore
}

\def\HEXAGONnoxy#1#2{
\save
=( 4, 0     )"A";
=( 2, 3.4641)"B";
=(-2, 3.4641)"C";
=(-4, 0     )"D";
=(-2,-3.4641)"E";
=( 2,-3.4641)"F";
\LEGS#2
\HEX#1
\restore    
}

\def\HEXAGON#1#2{
\begin{xy}
\HEXAGONnoxy{#1}{#2}
\end{xy}
}
\def\glue{\!\!\supset}
\def\path{\curvearrowleft}

\title{Phylogenetic toric varieties on graphs}
\author{Weronika Buczy\'nska
\thanks {This is shortened version of author's PhD thesis.  
Author was partially supported by the Texas Advanced Research
under Grant No. 010366-0054-2007. }
}

\begin{document}
\maketitle
\begin{abstract}
We define phylogenetic projective toric model of a trivalent graph as a generalization 
of a binary symmetric model of a trivalent phylogenetic tree.
Generators of the projective coordinate ring of the models of graphs with one cycle
are explicitly described.
The phylogenetic models of graphs with the same topological invariants
are deformation equivalent and share the same Hilbert function.
We also provide an algorithm to compute the Hilbert function.
\end{abstract}
\medskip
\footnotesize{
\noindent\textit{2010 Mathematics Subject Classification.}
Primary: 14M25; Secondary: 60J20, 92D15.\\
\noindent\textit{Keywords:} binary symmetric model, GIT quotient, Hilbert function.\\
}

\tableofcontents

\section{Introduction and background}\label{section_introduction}

The inspiration for this work are toric varieties arising 
in computational biology, 
or more precisely in phylogenetic algebraic geometry. 
The references to the subject include 
\cite{felsenstein_Inferring_Phylogenies},
\cite{pachter_sturmfels_alg_stat_for_comp_biology} and
\cite{semple_steel_phylogenetics}.

Markov models on phylogenetic trees are statistical models 
describing evolution.
They are usually defined as a subset of the probability simplex,
parametrized by a subset of matrices depending on the model.
Among them there are group-based models on phylogenetic trees.
These are special, as their projective versions, 
that is the Zariski closure of the parametrization in the 
complex projective space, are projective toric varieties.

We are interested in the simplest group-based models 
--- binary symmetric models, 
also called the Jukes-Cantor models, on trivalent trees. 
The object of study is the generalization 
of those models to trivalent graphs.

The article is organized in the following way:
in Section~\ref{section_motivation},
we give motivation to the subject
and point out references where the main object of our study appears.
Then we give a brief introduction to our main tools:  
we recall geometric invariant theory
in Section~\ref{section_GIT},
next, in Section~\ref{section_toric_in_wps}
we set the notation for projective toric varieties 
and in Section~\ref{section_quotient_by_subtorus}
provide the description of GIT quotient of
a projective toric variety by subtorus of its big torus.

In Section~\ref{section_graphs} we state 
combinatorial relations between the topological invariants 
of a trivalent graph and we prove that graphs with the same invariants 
are mutation-equivalent.
In Section~\ref{section_model}
we define the model of a trivalent graph
as a GIT quotient of a product of $\PP^3$'s 
indexed by inner vertices of the graph. 
This implies that the model is toric and comes with an embedding
into a weighted projective space. Our first result is Theorem~\ref{cone_gens} 
that lists the set of minimal generators of the 
projective coordinate ring of the model, 
when the underlying graph has the first Betti number at most one.
Section~\ref{section_defomation} contains our second result
--- models of mutation-equivalent graphs are deformation equivalent.

\begin{manualnothm}{\ref{deformation_equivalent}}
Geometric models of connected trivalent graphs with $n$  leaves and 
the first Betti number $g$  are deformation equivalent in 
the projective toric variety $\PP_{g,n}$, which is a quotient of 
$\PP^{2^{n+2g-1}-1}$ by a $g$-dimensional torus.
\end{manualnothm}

In Section~\ref{section_Hilber_function} we prove that the Hilbert functions of 
mutation-equivalent models are equal (Theorem~\ref{same_hs}) and finally
we compute these Hilbert functions explicitly.

\subsection {Motivation --- Markov models on phylogenetic trees}
\label{section_motivation}

A \emph{phylogenetic tree} is an acyclic connected graph 
with additional data attached to its edges and vertices. 
At a vertex $v$ there is a finite ordered set $A_v$ called an alphabet.
At an edge with ends $v$ and $w$ there is a doubly-stochastic matrix 
(all rows and columns sums are 1)
with the $(i,j)^{th}$ entry indicating the probability 
of the $i^{th}$ letter of $A_v$ being changed to the $j^{th}$ letter $A_w$.
To construct a \emph{Markov model} on a phylogenetic tree 
we first need to indicate a set of observable vertices, 
for example the leaves of the tree. Then the model 
is the subvariety of the probability simplex, 
parametrized by a subset of matrices that we only allow,
given by probabilities of observing letters at the observable vertices.
We consider symmetric models, which means we allow symmetric matrices. 
Typically the observable vertices are the leaves of the tree.

Apart form this real variety, one can consider its complex algebraic relaxation.
That is, the parameters are allowed to vary in a complex projective space
and we take the Zariski closure of the image. Then the model becomes 
a complex projective variety and can be studied by means of algebraic geometry.
Binary symmetric models have additional structure --- they are equipped 
with an action of a torus of dimension equal to the dimension of the model 
and thus they are projective toric varieties.
This is an especially nice class of varieties, which have 
a combinatorial description by lattice polytopes.
The geometry of the simplest group-based models --- binary symmetric models
with the restriction that the underlying tree is trivalent 
was the object of study of~\cite{buczynska_wisniewski}. 
In that paper  we described the corresponding 
lattice polytope and interpreted the models as 
a certain quotient of a product of three-dimensional projective spaces.

\subsection{Toric algebras of our graph models in the literature}
\label{section_literature}

We generalize the quotient description of the tree models 
introduced in~\cite{buczynska_wisniewski} and we again 
have a toric projective model, which this time is embedded 
in a weighted projective space. 
Such an embedding is always given by a graded lattice cone.
We denote the cone for a graph $\ccG$ by  $\tau(\ccG)$.

The way we associate a lattice cone to a trivalent graph 
appears also in the work of Manon~\cite{manon_conformal_blocks}.
He constructs a sheaf of algebras over the moduli stack ${\cal  M}_{g,n}$ 
of genus $g$ curves with $n$ marked points and our semigroup algebras $\CC[\tau(\ccG)]$
are obtained by some initial term deformations from algebras above the
most special points of  ${\cal  M}_{g,n}$ in Manon's construction.

Another place where our cones $\tau(\ccG)$ appear 
is the Jeffrey and Weitsmann's~\cite{jeffrey_weitsman} study 
of flat $SU(2)$-connections on a genus $g$ Riemann surface. 
In their context the trivalent graph $\ccG$ describes the geometry 
of the compact surface $\Sigma^g$ of genus $g$ and thus has no leaves. 
A subset of $\ZZ$-labellings of the graph, 
which are exactly points of our cone $\tau(\ccG)$, 
are in 1--1 correspondence with the number of Bohr-Sommerfeld fibers 
which is the central object of study in \cite{jeffrey_weitsman}. 
By the Verlinde formula, the number of those fibers equals 
the dimension of holomorphic sections of powers of
a natural line bundle on the moduli space of flat $SU(2)$ 
connections on $\Sigma^g$.
This number is a value of the Hilbert function of the toric model 
of a connected graph with  no leaves and the first Betti number $g$.

By Theorem~\ref{same_hs}, we know that the Hilbert function 
only depends on the topological invariants of the graph.

Although the model depends on the shape of the underlying trivalent tree,
once we restrict ourselves to trees with fixed number of leaves,
models of all of them are in the same irreducible component 
of the Hilbert scheme of projective varieties with fixed Hilbert polynomial.
This was proved by Sturmfels and Xu in \cite{sturmfels_xu}. 

Any trivalent graph is made by gluing together tripods, 
that is graphs $\TRIPOD$ with four vertices and three edges
attached to the central vertex. 
To construct the toric model we assign to every inner vertex
a copy of a three-dimensional complex projective space and
to every edge we assign an action of the one-dimensional complex torus $\CC^*$
on the product of all those $\PP^3$,
which corresponds to gluing two tripods along that edge.
The model $X(\ccG)$ of the trivalent graph $\ccG$ is 
a geometric invariant theory (GIT) quotient of product of the $\PP^3$ by 
the torus defined as a product of the $\CC^*$'s corresponding to 
the internal edges.
We also translate this description into language of projective toric varieties,
by writing the model $X(\ccG)$ as the projective spectrum 
of a semigroup ring $\CC[\tau(\ccG)]$. The underlying semigroup $\tau(\ccG)$
has a clear description in terms of the graph $\ccG$.

Three results of this article generalize our earlier results  
obtained in~\cite{buczynska_wisniewski} about
binary symmetric models of trivalent trees to phylogenetic graph models. 
First we describe the minimal $\ZZ$-generators of 
the semigroup $\tau(\ccG)$ when the graph $\ccG$ has 
the first Betti number at most one.
We also prove  that models of graphs with the same 
discrete invariants are deformation equivalent
and lastly that they share the same Hilbert function.

\subsection{Geometric invariant theory.}\label{section_GIT}
We use geometric invariant theory  for a normal projective variety $X$ 
with an action of an algebraic torus $\TT$. 
Our main reference is Section~5~and~6 of~\cite{bialynicki_enc},
although the setup we use may seem to be slightly more general 
then the  one found in~\cite{bialynicki_enc}. 
This is because instead of linearizing only with respect to a line bundle 
we allow ample Weil divisors.
We will explain the necessary modifications and show how this 
does not affect the basic theory.

\begin{defin}
A divisor $\ccL$ is an \textbf{ample Weil divisor} if 
some positive multiple $n\ccL$ is an ample line bundle.
\end{defin}
Given an ample Weil divisor $\ccL$ we have the ring
\[
R(X,\ccL):=\bigoplus_{p=0}^\infty H^0(X,\ccO(p\ccL)),
\]
which is the projective coordinate ring of $X$ embedded 
into a weighted projective space by the linear system $|\ccL|$. 
This is completely analogous, see \cite{reid_graded_rings_birational_geometry},
to the standard way of describing 
embedding of $X$ into a projective space in when $\ccL$ 
is a very ample line bundle, see \cite[Section II.2]{hartshorne}.
We discuss those facts in~\ref{section_toric_in_wps}.

We denote by
\[
R_p(X,\ccL):=R(X,p\ccL)
\]
the ring given by a multiple of $\ccL$.
When the divisor $\ccL$ is clear,
we write $R(X)$ and $R_p(X)$ an instead of $R(X,\ccL)$ and $R(X,p\ccL)$.
\begin{lem}[Veronese embedding]\label{veronese}
Let $X=\Proj R(X,\ccL)$ be a projective variety
with a distinguished ample Weil divisor $\ccL$.
Then for any positive integer $p$ 
the inclusion $R_p(X,\ccL)\mono R(X,\ccL)$ induces an isomorphism
$$\Proj R(X,\ccL) \simeq \Proj R_p(X,\ccL)$$
\end{lem}
\begin{prf}
See \cite[Chapter 2, Theorem~2.4.7]{EGAII},
and also \cite[Ex. II 5.13]{hartshorne}.
\end{prf}

We briefly recall the definitions of quotients simplifying 
the statements from \cite{bialynicki_enc}. 
Let $G$ be a reductive algebraic group acting on a variety $X$.
A $G$-equivariant morphism  $\pi : X \ra Y$ is called
a \textbf{semi-geometric quotient} if
\begin{itemize}
\item the image of every closed orbit is closed,
and this property is invariant under base change
\item $\pi$ is surjective and images of disjoint, closed orbits are disjoint
and this property is invariant under base change
\item $\pi_*(\ccO^G_X)=\ccO_Y$.
\end{itemize}
The map $\pi : X \ra Y$ is called a \textbf{good quotient}
if it is both affine and semi-geometric. 
A map $\pi : X \ra Y$ is affine if preimage of any affine subset of $Y$ 
is affine.

This definition implies that the good quotient is a categorical quotient, 
and thus unique \cite[Rmk.3.1 and Thm~3.2]{bialynicki_enc}.

To construct the quotient, we need some more definitions.
\begin{defin}
Let $G$ be a reductive algebraic group acting on $X$.
An ample Weil divisor $\ccL$ on $X$ is \textbf{G-linearized} 
if the action of $G$ can be lifted to 
the projective coordinate ring $R(X,\ccL)$.
More precisely, the action of $G$ on $X=\Proj R(X,\ccL)$ 
induced by the lift, is the action of $G$ on $X$.
\end{defin}

When $\ccL$ is a line bundle, a $G$-linearization
is an action of $G$ on $\ccL$ which is linear on fibers
and agrees with its action on $X$. 
We recall the definition of 
a $\ccL$-semistable point~\cite[\S~6.1]{bialynicki_enc},
which is independent of the multiple of $\ccL$.

\begin{defin}
Let $\ccL$ be an ample Weil divisor.
A point $x \in X$ is $\ccL$-\textbf{semistable} if there exists 
a $G$-invariant section of a positive multiple of $\ccL$ 
that does not vanish at $x$.
We denote by $X^{ss}$ the set of all semistable points.
The \textbf{GIT quotient} of $X$ by $G$ with respect 
to the linearization $\ccL$~is
\[
X \git G = X^{ss} \git G,
\] 
where $X^{ss} \git G$ denotes the good quotient of the set of
semistable points with respect to $\ccL$ by the action of $G$.
The rational map $X \ra X\git G$ is called GIT quotient map and 
the regular map $X^{ss} \ra X \git G$ is a good quotient map.
\end{defin}

\begin{rmk}
Both \cite{mumford_fogarty_kirwan_geometric_invariant_theory}
and \cite{bialynicki_enc}  assume that the section in the above definition
has an affine support, but as we consider only
ample $\ccL$ all its sections automatically have affine support.
\end{rmk}

The following theorem says that in the above situation a GIT quotient exists 
and is the projective spectrum of the ring of invariants. 
\begin{thm}\label{GIT}
Let $G$ be a reductive group acting on projective varieties $X$ and $Y$.
\begin{enumerate}
\item\label{quotient}
Let $\ccL$ be a $G$-linearized ample  Weil divisor on $X$.
Then \[X  \git  G = \Proj R(X,\ccL)^G.\]
\item\label{subquotient}
Let $X \mono Y$ be a $G$-equivariant embedding,
with $Y=\Proj R$ and $X=\Proj R / I$, 
where $I$ is the homogeneous ideal of $X$ in $Y$.
Then $I$ is generated by invariants $f_1, \ldots, f_j$,
the map $X \git G \mono Y \git G$ is an embedding and
the homogeneous ideal of $X \git G$ in $R^G$ 
is also generated by $f_1, \ldots, f_j$.
\end{enumerate}
\end{thm}

\begin{prf}
In \cite[\S5]{bialynicki_enc} the affine quotient is defined 
for an affine variety $X$ as a spectrum of invariants and 
\cite[Thm~5.4]{bialynicki_enc} says it is a good quotient.
Then in \cite[\S7]{bialynicki_enc} there is a characterization of 
a good quotient $\pi : X \ra Y$ as a locally affine quotient, that is 
the map $\pi$ is a good quotient if and only if
for any open affine subset $U \subset Y$ the restriction 
of the map $\pi$ to $\pi^{-1}(U)$ is an affine quotient.

We know from \cite[Thm~6.2.1]{bialynicki_enc} that~\ref{quotient}
is true for projective space and $\ccL=\ccO(1)$. 
Additionally for an arbitrary variety $X=\Proj R$ and 
a very ample line bundle $\ccL$ 
that defines an equivariant embedding into $\PP^n$ 
by \cite[\S6.3]{bialynicki_enc} we have:
\begin{itemize}
\item $X^{ss}={(\PP^n)}^{ss} \cap X$
\item the restriction of the good quotient morphism 
$\pi \colon (\PP^n)^{ss} \ra (\PP^n)^{ss} \git G$
is the good quotient morphism on~$X$.
\end{itemize}

We choose an affine covering $\cal U$ of the quotient $\PP^n \git G$
by sets of the form $\{f \ne 0\}$, where $f$ is $G$-invariant.
Then $\{\pi^{-1}(U) | U \in \cal U \}$ is 
an affine $G$-invariant covering of $\PP^n$, 
since $\pi$ is an affine map. 
By \cite[Thm~5.3]{bialynicki_enc} and 
the Hilbert-Nagata theorem \cite[Thm~5.2]{bialynicki_enc}
for each  $U \in \cal U$ we know that
$(\pi^{-1}(U) \cap X) \git G$ is a spectrum of invariants
\[
(\pi^{-1}(U) \cap X) \git G = \Spec (R[f^{-1}]^0)^G
\]
Also each $\pi^{-1}(U) \cap X$ has form $\{g=\pi^*f|_X \neq 0\}$.
To see that $X \git G = \Proj R^G$, we only need to see that
$\left( R^G[g^{-1}] \right)^0 = (R[g^{-1}]^0)^G $, which
is true since $g$ is $G$-invariant.
Thus~\ref{quotient} holds for $X$ projective and $\ccL$ very ample.

When $\ccL$ is an ample Weil divisor, we use Lemma~\ref{veronese}
to replace $\ccL$ by its $k^{th}$ power. More precisely,
we choose $k$ such that $k\ccL$ is a very ample line bundle,
and since the set of semistable points does not depend on the 
multiple of $\ccL$ we have
\[
X \git G = \Proj R_p(X,\ccL)^G = \Proj(R(X,\ccL)^G)_p = \Proj R(X,\ccL)^G
\]
where the first equality holds, because $k\ccL$ is very ample,
the second holds because the action of $G$ preserves gradations, 
and in the third we use Lemma~\ref{veronese}.

For the proof of~\ref{subquotient} by Hilbert-Nagata theorem, 
we know that $I$ is generated by invariants.
Since a good quotient is a locally affine quotient 
the map $X \git G \ra Y \git G$ is an embedding.
The last part of~\ref{subquotient}
follows from~\cite[Thm.~5.2]{bialynicki_enc}.
\end{prf}

\subsection{Toric varieties  in weighted projective spaces.}
\label{section_toric_in_wps}

The embedding of a projective toric variety $X$ into a projective space 
is described by a polytope $\Delta$ with integral vertices. 
If we scale the polytope by an integer $k$, then we do not change the variety.
The resulting embedding changes by composing it with $k^{th}$ Veronese
embedding. The sum of all positive multiples forms a semigroup 
(or a graded cone).
In this situation $X = \Proj \bigoplus_{k\in \NN} \CC [k \Delta]$.
When the ambient space is a weighted projective space the embedding 
is given by a graded cone with a set of (minimal) $\ZZ$-generators, 
which are no longer in the first degree.
Every section of this cone determined by the grading is a rational polytope.

\begin{defin}
A \textbf{weighted projective space} $\PP(a_0,\ldots,a_n)$ 
with weights $(a_0,\ldots,a_n)$ where each $a_i \in \NN$ 
is a positive integer is the GIT quotient of 
the affine space $\CC^{n+1}$ by the action 
\[
t \cdot (x_0,\ldots,x_n) = (t^{a_0}\cdot x_0,\ldots,t^{a_n} \cdot x_n)
\]  
Again, the non-stable locus is the point $0 \in \CC^{n+1}$.
\end{defin}
We can assume that the  greatest common divisor
of the weights is one --- this is by substitution
$t \mapsto t^{\gcd(a_0,\ldots,a_n)}$.
Let us pick integers $b_0,\ldots,b_n$ such that 
$\sum_{i=0}^n b_i \cdot a_i = 1$.
The weighted projective space has the sheaf $\ccO_{\PP}(1)$, 
which corresponds to the ample Weil divisor  $\sum_{i=0}^n b_i \cdot (x_i)$,
where the  $(x_i)$'s are the divisors corresponding to the coordinates.  

Weighted projective spaces are often singular.
They have quotient singularities coming from finite abelian group actions.

\begin{defin}\label{graded_lattice}
  Given a lattice $M$, we can associate with it in a non-unique way 
a \textbf{graded lattice}  $\Mg$, which equipped with degree map
$\deg \colon \Mg \ra \ZZ$ --- the projection to the first coordinate
and fits into the exact sequence
\[
\xymatrix{
0 \ar@{->}[r] &
M \ar@{->}[r] &
\Mg \ar@{->}[r]^{\deg}  &
\ZZ \ar@{->}[r]  \ar@<+.6ex>[l]^{s} &
0\\
}\]
We also fix a splitting $s : \ZZ \ra \Mg$ of the exact sequence 
which is equivalent to a choice of the $0$ element 
in the $M \simeq (1,M) \subset \Mg$. 
The choice of the splitting $s$ corresponds to a choice of linearization 
of the action of the torus $\Spec (\CC[M^{\vee})$ on itself.
\end{defin}

\begin{defin}
  A \textbf{graded lattice cone} $\tau$ is 
a rational, convex, polyhedral cone in a graded lattice $\Mg$, 
with all elements having non-negative degree: $\deg(\tau) \subset \NN$,
and the zero gradation consists of one element: $\deg^{-1}(0) \cap \tau={ 0}$.
Convex here means that 
$\tau \otimes_{\ZZ} \RR_{+} \subset \Mg \otimes_{\ZZ} \RR$ is convex.
Equivalently, $\tau$ is a (saturated) sub-semigroup of 
the free abelian group $\Mg$ with finite set of $\ZZ$-generators 
all having positive degrees
where  the neutral element is the only one of degree zero.
\end{defin}

\begin{prop}\label{toric_in_wps}
A toric variety $X\subset\PP(a_0,\ldots,a_n)$ is described by 
its fan and an  ample Weil divisor $\ccO_X(1)$ or 
equivalently by an isomorphism $X \simeq \Proj \CC[\tau]$ 
where $\tau$ is a graded, rational, convex, polyhedral cone
in a graded lattice $\Mg$. Then $\CC[\tau]=R/I$, where 
$R$ is the homogeneous coordinate ring of $\PP(a_0,\ldots,a_n)$
and $I$ is the homogeneous ideal of $X$.
\end{prop}
\begin{prf}[Idea of the proof]
The correspondence between the ample Weil divisor and 
the graded cone is the following. 
The degree $k$ sections of the sheaf associated to the Weil divisor
form the $k^{th}$ section of the cone $\tau$. 
To go the other way, we have an isomorphism $X \simeq \Proj \CC[\tau]$ 
and then the  ample Weil divisor is $\ccO_{X}(1)$ ---
the pull-back of $\ccO_{\PP(a_0,\ldots,a_n)}(1)$
from the ambient weighted projective space $\PP(a_0,\ldots,a_n)$.  
\end{prf}

The choice of the $0$ element of $M \simeq (1,M) \subset \Mg$ 
in Definition~\ref{graded_lattice} is a choice 
of linearization of the action of the torus of $X$, 
which extends the action of the torus on itself.

\begin{defin}
  A \textbf{graded product} $\Mg_1 \times_{g} \Mg_2$ of the graded lattices
$\Mg_1$ and $\Mg_2$ 
is the fiber product over their degree maps, or equivalently
the hyperplane $\deg_1=\deg_2$ in the product $\Mg_1 \times \Mg_2$:
\[
\xymatrix{
\Mg_1 \times_g \Mg_2 \ar@{->}[r] \ar@{->}[d] & 
\Mg_2 \ar@{->}[d]^{\deg_2}\\
\Mg_1 \ar@{->}[r]^{\deg_1} & \ZZ \\
}\]
\end{defin}

\begin{defin}
A \textbf{graded product} $\tau_1 \times_g \tau_2$ of graded cones 
$\tau_1$ and $\tau_2$ is the fiber product over their degree map,
or equivalently intersection of the hyperplane $\deg_1 = \deg_2$ with
the product cone $\tau_1 \times \tau_2$:
\[
\xymatrix{
\tau_1 \times_g \tau_2 \ar@{->}[r] \ar@{->}[d] & 
\tau_2 \ar@{->}[d]^{\deg_2}\\
\tau_1 \ar@{->}[r]^{\deg_1} & \ZZ \\
}\]
\end{defin}

\begin{lem}\label{product_of_varieties_is_a_graded_product_of_cones}
If $\tau_1 \subset \Mg_1$ and $\tau_2 \subset \Mg_2$ are 
graded rational convex polyhedral cones, and
\[
X_1 = \Proj\CC[\tau_1] \qquad X_2 = \Proj \CC[\tau_2]
\]
are corresponding projective toric varieties, 
then the product of these varieties under the Segre embedding 
\[
X_1 \times X_2 = \Proj \CC[ \tau_1 \times_g \tau_2]
\] 
corresponds to the graded product of the cones $\tau_1$ and $\tau_2$.
\end{lem}
\begin{prf}
By definition  $X_1 \times X_2$ under Segre embedding is equal to
\[
\Proj \bigoplus_{i=0}^{\infty} H^0\bigl(
X_1 \times X_2, \ccO_{X_1}(i) \boxtimes \ccO_{X_2}(i)
\big).
\]
For each $i$ we know that the sections of 
this exterior tensor product are spanned by 
the product of the $i^{th}$ graded pieces of the cones
\[
H^0\bigl(
X_1 \times X_2, \ccO_{X_1}(i) \boxtimes \ccO_{X_2}(i)
\big) =
\CC \cdot \Bigl( 
\bigl(\tau_1 \cap \deg^{-1}(i)\bigr) \times 
\bigl(\tau_1 \cap \deg^{-1}(i)\bigr)
\Bigr)
\]
We conclude the lemma by summing the above equality over all $i$'s to get
\[
\bigoplus_{i=0}^{\infty} H^0\bigl(
X_1 \times X_2, \ccO_{X_1}(i) \boxtimes \ccO_{X_2}(i)
\big)
= \CC [\tau_1 \times_g \tau_2]
\]
as required.
\end{prf}

\subsection{Quotient of a projective toric variety by subtorus.}
\label{section_quotient_by_subtorus}

We show that the quotient of a projective toric variety 
$X=\Proj \CC[\tau]$ by a subtorus is described by 
an appropriate linear section of $\tau$.

We first recall some facts about toric varieties from \cite{fulton}.
If $X=\Proj \CC[\tau]$, where $\tau \subset \Mg$ is a graded cone,
then the torus of $X$ is given by an isomorphism 
$\TT=\Hom(M, \ZZ)\otimes_{\ZZ}\CC^*$.
The lattice $N=\Hom(M, \ZZ)$ is the lattice of
one-parameter subgroups of $\TT$.
For any subtorus $\TT'$ of the torus $\TT$ there are
corresponding maps of lattices:
the projection $M \epi M'$ of the monomial lattices
and the embedding $N' \mono N$ of the lattice of one-parameter subgroups.

A linearization of the action of the torus $\TT$ on $X$, 
induces a linearization of the action of any subtorus $\TT' \subset \TT$, 
by restricting the action. In this setting Theorem~\ref{GIT} implies 
the following result.

\begin{thm}\label{toricGIT}
Let $\tau$ be a graded cone in a lattice $\Mg$ and $X=\Proj \CC[\tau]$ 
the corresponding toric variety. 
Let  $j:\TT' \mono \TT$  be a subtorus of the torus $\TT$, 
with $j_* : N' \mono N$, $j^* : M \epi M'$ 
and $\id \oplus j^* : \Mg \ra (M')^{gr}$ the corresponding lattice maps.
Then there exists a good quotient and it is equal to
\[
X \git \TT' = \Proj \CC[\tau]^{\TT'},
\]
where
\[
\CC[\tau]^{T'} = \CC[\tau \cap (\ZZ\oplus \ker(j^* : M\epi M'))]
=\CC[\tau \cap (\ZZ \oplus N'^{\perp})].
\]
Moreover the quotient $X \git \TT'$ is polarized by 
$\ccO_{X\git \TT'}(1)$ in a natural way.
\end{thm}

The following example shows that we do need ample Weil divisors,
not only ample line bundles.

\begin{ex}
Let $X$ be the good GIT-quotient
\[
\pi : \PP^3 \times \PP^3 \ra X = \PP^3 \times \PP^3) \git (\CC^*)^3  
\]
of the product of two projective three-spaces by an action of
three-dimensional torus acting with weights
\[
\begin{bmatrix*}[r]
(0& 1& 1& 0) &\times& (0& -1& -1& 0) \\
(0& 1&-1& 0) &\times& (0& 0& 0& 0) \\
(0& 0& 0& 0) &\times& (0& 1&-1& 0)
\end{bmatrix*}
\]
linearized with respect to the line bundle 
$\ccO_{\PP^3}(1) \boxtimes \ccO_{\PP^3}(1) $.
We will see later that $X=X \left( \DUMBBELLfont \right)$
is the model of the trivalent graph $\DUMBBELLfont$ 
with the first Betti number two with three edges
and is a projective toric variety by Theorem~\ref{toricGIT}.
The sheaf $\ccO_X(1)=\pi_*(\ccO_{\PP^3}(1) \boxtimes \ccO_{\PP^3}(1))$ 
is not a locally free \mbox{$\ccO_X$-module}
because the associated divisor is not Cartier. 
To verify it we can use a computer algebra system, 
for example magma \cite{magma} as follows.
Since any divisor on a toric variety is linearly equivalent 
to a $\TT$-invariant divisor, we identify a divisor with a corresponding 
$\ZZ$-combination of primitive elements of the rays of the fan.
Thus we only need to check if the $\ZZ$-combination corresponding 
to $\ccO_X(1)$  yields a piecewise linear function on the fan, 
which  by~\cite{oda} is equivalent to our 
$\TT$-invariant Weil divisor being Cartier.
\end{ex}

\section{Phylogenetic  models on trivalent graphs}%
\label{section_model}
\subsection{Trivalent graphs.}\label{section_graphs}

We define topological invariants of trivalent graphs and show 
any two graphs with the same invariants are equivalent by applying 
appropriate mutations, which we introduced in \cite{buczynska_wisniewski}.
We \emph{do not} assume that our graphs are connected.

\begin{defin}
A \textbf{graph} $\ccG$ is set $\ccV$ of vertices and set $\ccE$ of edges  
together with the unordered boundary map $\d: \ccE \ra \ccV^{\otimes 2}$, 
where $\ccV^{\otimes 2}$ is the set of unordered pairs of vertices.
We write $\d(e) = \{ \partial_1(e), \partial_2(e)\}$ and say that $v$ is
an \textbf{end of the edge} $e$ if $v \in \d(e)$.
A vertex incident to exactly one edge is a \textbf{leaf}.
The \textbf{set of leaves} is denoted by $\ccL$ and 
the \textbf{number of leaves} by $n$. 
If a  vertex is  not a leaf, it is called an \textbf{inner vertex}. 
An edge incident to a leaf is a \textbf{petiole} and 
$\ccP$ is the set of petioles.
We write $\comp\ccG$ for the set of \textbf{connected components}
of the graph and $|\comp\ccG|$ for the number of components.
We denote by $g$ for the \textbf{first Betti numbers} of graph, which is 
the rank of the first group homology of the graph viewed as a CW-complex.
A graph is \textbf{trivalent} if every inner vertex has valency three.
Valency of a vertex $v$ is the number of connected components 
of a sufficiently small neighborhood of $v$ with $v$ removed.
A trivalent graph with no cycles is a \textbf{trivalent tree}.
\end{defin}

When discussing more then one graph
instead of $\ccV$, $\ccE$, $n$,...
we will write $\ccV(\ccG)$, $\ccE(\ccG)$, $n(\ccG)$,~etc.

\begin{rmk}
Our graphs are \emph{not oriented}, nevertheless we write $\partial_1(e)$
and $\partial_2(e)$ for the vertices adjacent to the edge $e$.
This makes it easier to talk about ``the other end of $e$''.
\end{rmk}

We call the unique trivalent tree $\TRIPOD$ 
with a single inner vertex and three leaves 
the \textbf{elementary tripod}. It has three edges $e_1, e_2, e_3$.
Any trivalent graph is built of elementary tripods in the following way:
given a trivalent graph $\ccG$ and any inner vertex $v \in \ccV(\ccG)$ 
we pick a copy of elementary tripod $\TRIPOD_v \simeq \TRIPOD$ 
and an map $i_v : \TRIPOD_v \ra \ccG$ 
which sends the central vertex of $\TRIPOD_v$ to $v$
and locally near $i^{-1}(v)$ is an embedding.
We present the graph $\ccG$ as a disjoint union of the graphs $\TRIPOD_v$ 
with appropriate identification of edges
\begin{equation}\label{graph_from_tripods}
\ccG = \bigsqcup_{v \in \ccV \setminus \ccL} \TRIPOD_v / 
\{  i_{\partial_1(e)}^{-1}(e)\sim i_{\partial_2(e)}^{-1}(e) 
\}_{e\in\ccE\setminus\ccP}
\end{equation}
This construction mirrors how the model of the graph is constructed, 
as we will see in Definition~\ref{defin_as_quotient}.
\begin{ex}
On the Figure~\ref{pic_building_graph}
we give an example of the above presentation of a trivalent graph
for a graph with the first Betti number one and two leaves.
\myVSPACEfigure\begin{center}
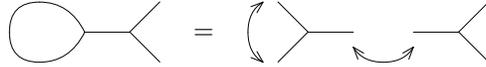
  \begin{minipage}{\textwidth}
\[
\begin{xy}
(0,0);(0,2):;
( 0, 0)*{} ; ( 2, 3)*{} **\crv{~**\dir{-} (2,1)};
( 2, 3)*{} ; ( 0, 5)*{} **\crv{~**\dir{-} (2,5)};
(-2, 3)*{} ; ( 0, 5)*{} **\crv{~**\dir{-} (-2,5)};
( 0, 0)*{} ; (-2, 3)*{} **\crv{~**\dir{-} (-2,1)};
( 0, 0)*{} ; ( 0,-3)*{} **\dir{-};
( 0,-3)*{} ; (-2,-5)*{} **\dir{-}; 
( 0,-3)*{} ; ( 2,-5)*{} **\dir{-};
\end{xy}\quad=\quad 
\begin{xy}
(0,0);(0,2):;       
( 2, 4)*{} ; ( 0, 5)*{} **\crv{~**\dir{-} (1.4,5)};
( 2, 4)*{} ; ( 1.7, 4.8)*{} **\dir{-};
( 2, 4)*{} ; ( 1.3, 4.4)*{} **\dir{-};
(0,0);(-1,0):(0,-1)::;
( 2, 4)*{} ; ( 0, 5)*{} **\crv{~**\dir{-} (1.4,5)};
( 2, 4)*{} ; ( 1.7, 4.8)*{} **\dir{-};
( 2, 4)*{} ; ( 1.3, 4.4)*{} **\dir{-};
(0,-2);(1,-2):;     
( 0, 0)*{} ; ( 0,-3)*{} **\dir{-};
( 0,-3)*{} ; (-2,-5)*{} **\dir{-}; 
( 0,-3)*{} ; ( 2,-5)*{} **\dir{-};
(-3,2);(-3,1):;       
( 2, 4)*{} ; ( 0, 5)*{} **\crv{~**\dir{-} (1.4,5)};
( 2, 4)*{} ; ( 1.7, 4.8)*{} **\dir{-};
( 2, 4)*{} ; ( 1.3, 4.4)*{} **\dir{-};
(0,0);(-1,0):(0,-1)::;
( 2, 4)*{} ; ( 0, 5)*{} **\crv{~**\dir{-} (1.4,5)};
( 2, 4)*{} ; ( 1.7, 4.8)*{} **\dir{-};
( 2, 4)*{} ; ( 1.3, 4.4)*{} **\dir{-};
(2,3);(2,2):;      
( 0, 0)*{} ; ( 0,-3)*{} **\dir{-};
( 0,-3)*{} ; (-2,-5)*{} **\dir{-}; 
( 0,-3)*{} ; ( 2,-5)*{} **\dir{-};
\end{xy}
\]
\captionof{figure}{Building a trivalent graph from tripods}
\label{pic_building_graph}
\end{minipage} \end{center}\myVSPACEfigure
\end{ex}

\begin{lem}\label{graph_invariants}
In any trivalent graph with $n$ leaves and first Betti number $g$
the following holds
\begin{enumerate}
\item $|\ccV|,|\ccE| \ge n$,
\item\label{trivalent_equation} $2|\ccE|=3|\ccV|-2n$,
\item\label{euler_characteristic} $|\ccV|-|\ccE|=|\comp\ccG|-g$.
\end{enumerate}
Thus, any three of the numbers $|\ccV|$, $|\ccE|$, $n$, $g$, $|\comp\ccG|$ 
determine the other two.
\end{lem}
\begin{prf}
To prove \ref{trivalent_equation} let us count pairs of consisting of 
a vertex and an adjacent edge. On one hand we will count every edge twice.
On the other hand every inner vertex has three incident edges so we have 
$3(|\ccV|-n)$ pairs and another $n$ pairs come from leaves 
which totals to $3|\ccV|-2n$.
Equation \ref{euler_characteristic} counts the Euler characteristic 
$|\comp\ccG|-g$ of the graph.
\end{prf}
Here we introduce operations of gluing two leaves of a graph, cutting an edge into two new edges
and taking a disjoint sum of two graphs.

\begin{defin}\label{definition_operations_on_graphs}
We will use the following three constructions of trivalent graphs.
\begin{itemize}

\item  $\ccG \sqcup \ccG'$ is the  \textbf{disjoint sum} of the given graphs  
$\ccG$ and $\ccG'$.

\item $\ccG^{l_1}_{l_2} \glue$ is the graph obtained from a given graph $\ccG$ 
with two distinguished leaves $l_1, l_2 \in \ccL(\ccG)$ 
by \textbf{gluing the two leaves} $l_1$ and $l_2$,
or more precisely by removing $l_1$ and $l_2$
and identifying the edge incident to $l_1$ with the edge incident to $l_2$.

\item $\ccG_l \star \ccG'_{l'}$ a \textbf{graft} of given graphs 
$\ccG$ and $\ccG'$ each with a distinguished leaf. 
Figure~\ref{pic_graft} is a schematic picture of this construction.

\myVSPACEfigure\begin{center}
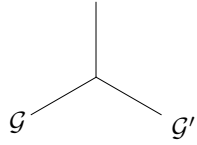
\begin{minipage}{\textwidth}
\[
\begin{xy}
(0,0)*{}="0";
(0,10)*{} ="A"; 
(-\fiverootthree,-5)*{}="B" ; 
( \fiverootthree,-5)*{}="C" ; 
(0,-4)="a"; ( \rootthree,1)="b"; (-3,  \rootthree)="c";
"0" ; "A" **@{-}-"a";
"0" ; "B" **@{-}-"b"*{\ccG};
"0" ; "C" **@{-}-"c"*{\ccG'};
\end{xy}
\]
\captionof{figure}{Graft of two graphs}
\label{pic_graft}
\end{minipage} \end{center}\myVSPACEfigure

The new graph can be written as 
\[
\ccG_l \star \ccG'_{l'} = 
\left( \left((\ccG \sqcup \ccG')^l_{e_1} \glue\right) \sqcup
\ccG' \right)^{l'}_{e_2} \glue
\]

\item $\ccG^e$ is  the graph  obtained from the given graph $\ccG$ by 
\textbf{cutting an internal edge} $e \in\ccE(\ccG)\setminus\ccL(\ccG)$.
More precisely we replace $e$ by two new edges $e_1$ and $e_2$ with
$\partial_1(e_1):=\partial_1(e)$ and $\partial_1(e_2):=\partial_2(e)$. 
There are two new leaves in $\ccG^e$, which are the free ends of the new edges
$\partial_2(e_1)$ and $\partial_2(e_2)$.
\end{itemize}
\end{defin}

\begin{defin}\label{definition_cycle_edge_pending_tree_polygon_graph}
An edge $e \in \ccE$ is called a \textbf{cycle edge} if 
it is not a petiole and removing it  does not disconnect 
the connected component of the graph that contains $e$. 
An edge $e \in \ccE$ is called \textbf{cycle leg} if 
it is incident to a cycle edge but 
it is not a cycle edge.
A vertex  $v \in \ccV$ is called \textbf{cycle vertex} if 
it is an end of a cycle edge.
We draw example of those on Figure~\ref{cycle_edge_vertex_leg}.

\myVSPACEfigure\begin{center}
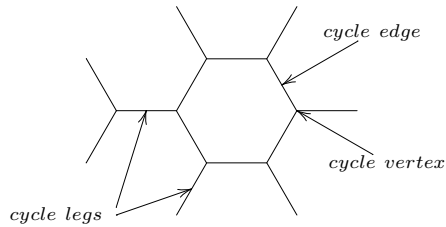
\begin{minipage}{\textwidth}
\[
\begin{xy}
0;<2mm,0mm>:;
\HEXAGONnoxy{------}{------}
\save
<0mm,0mm>;\POS<1.5mm,0mm>:; 
\POS"D"+"D"="DD"; \POS "E"+"E"="EE";
<0mm,0mm>;<1mm,0mm>:; \POS"A"+"B"="AB";
\restore\save
0;<2mm,0mm>:;
\restore\save
{\ar "AB"+"AB"+"AB"*{\scriptstyle{cycle\ edge}} ; "AB"*{}} 
\restore \save
{\ar "A" +"A" + "F"*{\scriptstyle{cycle\ vertex}} ; "A"*{} }
\restore \save
{\POS{"D"+"D" + "E" + "E"}*{\scriptstyle{cycle\ legs}}}
\restore \save
\drop{\ar "D" + "E" + "E" + (0,0.5); "DD"*{}}
\restore \save
\drop{\ar "D" +"E" + "E" ; "EE"*{}}
\restore \save
{\ar@{-} "D"+"D"+"C" ; "D"+"D" }
\restore \save
{\ar@{-} "D"+"D"+"E" ; "D"+"D" }
\end{xy}
\] 
\captionof{figure}{Cycle edge, cycle leg and cycle vertex of a graph}
\label{cycle_edge_vertex_leg}
\end{minipage} \end{center}\myVSPACEfigure
A \textbf{path} is a sequence of distinct edges $e_0,\ldots,e_m$
with  $\partial_2(e_i)=\partial_1(e_{i+1})$ 
for all $i \in \{0,\ldots, m-1\}$, moreover
$\partial_1(e_0)$ and $\partial_2(e_m)$ are both leaves or they are either equal.
In the latter case, the path is called a \textbf{cycle}.
Paths are \textbf{disjoint} if they have no common vertices.
A \textbf{network} is a union of disjoint paths.
For consistency we say that the empty set is also a network.
A \textbf{cycle} is a minimal sequence of cycle edges.
A cycle of length one is a \textbf{loop}.
On Figure~\ref{pic_path_cycle_network_loop} we draw examples for 
each those sequences.

\myVSPACEfigure\begin{center}
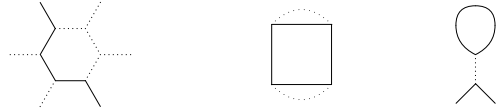
  \begin{minipage}{\textwidth}
\[
\xymatrix{
\begin{xy}
\HEXAGONnoxy{..---.}{..-..-} 
\end{xy}
\qquad\qquad\qquad
\begin{xy}
=( 1, 1)"A";
=(-1, 1)"B";
=(-1,-1)"C";
=( 1,-1)"D";
(0,0);(4,0):;
"A";"B" **@{-};
"B";"C" **@{-};
"C";"D" **@{-};
"D";"A" **@{-};
"A";"B" **\crv{~**\dir{.} (0, 2)};
"C";"D" **\crv{~**\dir{.} (0,-2)};
\end{xy}
\qquad \qquad
\begin{xy}
(0,0);(1.3,0):;
\LM-.--[][]
\end{xy}
}
\]
\captionof{figure}{A path, a cycle and a network containing a loop}
\label{pic_path_cycle_network_loop}
\end{minipage} \end{center}\myVSPACEfigure
A graph $\ccG$ is called a \textbf{polygon graph} if it has $2k$ edges 
of which $k$ form the only cycle of $\ccG$ and the remaining 
$k$ edges are cycle legs.
If $\ccG$ is any graph,  $e \in \ccV(\ccG)$ a non-cycle edge and 
after cutting $e$ we get a decomposition 
$\ccG^e=\ccG_0 \sqcup \ccG_1$ where $\ccG_1$ is a tree then we call
$\ccG_1$ a \textbf{pendant tree}.
\end{defin}

Figure~\ref{pic_mutation} shows the three trivalent trees 
with one internal edge $e$ and four labeled leaves.

\myVSPACEfigure\begin{center}
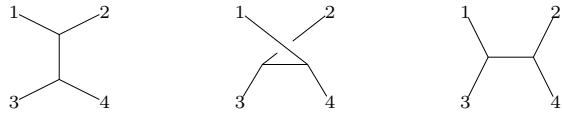
  \begin{minipage}{\textwidth}
\[
\begin{xy}
=(-3  , 3  )"A";
=( 3  , 3  )"B";
=(-3  ,-3  )"C";
=( 3  ,-3  )"D";
=( 0  , 1.5)"Q1";
=( 0  ,-1.5)"Q2";
=(-1.5, 0  )"O1";
=( 1.5, 0  )"O2";
=( 1.5, -0.5)"P1";
=(-1.5,-0.5)"P2";
=(-0.1, 0)"left";
=( 0.1, 0)"right";
0;<2mm,0mm>:;
"A"*{\scriptstyle 1} ; "Q1" **@{-};
"B"*{\scriptstyle 2} ; "Q1" **@{-};
"C"*{\scriptstyle 3} ; "Q2" **@{-};
"D"*{\scriptstyle 4} ; "Q2" **@{-};
"Q1"    ; "Q2" **@{-};
<30mm,0mm>;<32mm,0mm>:;
"A"*{\scriptstyle 1} ; "P1" **@{-}; {\ar@{-} 
"B"*{\scriptstyle 2} ; "P2" 
|!{"A";"P1"}\hole};
"C"*{\scriptstyle 3} ; "P2" **@{-};
"D"*{\scriptstyle 4} ; "P1" **@{-};
"P1"; "P2" **@{-};
<60mm,0mm>;<62mm,0mm>:;
"A"*{\scriptstyle 1} ; "O1" **@{-};
"B"*{\scriptstyle 2} ; "O2" **@{-};
"C"*{\scriptstyle 3} ; "O1" **@{-};
"D"*{\scriptstyle 4} ; "O2" **@{-};
"O1"; "O2" **@{-};
\end{xy}
\]
\captionof{figure}{The three trees with  four labeled leaves}
\label{pic_mutation} 
 \end{minipage} \end{center}\myVSPACEfigure

Let $\ccG$ be a trivalent graph and let $e$ be an internal edge
which is not a loop (the ends of $e$ are not identified).
Then a neighborhood of $e$ in $\ccG$ is a trivalent tree with four leaves.
\begin{defin}
A \textbf{mutation} of a trivalent graph $\ccG$ \textbf{along an edge $e$} 
is a graph $\ccG'$ which is obtained from $\ccG$ 
by removing a neighborhood of $e$ (which is a trivalent tree with four leaves) 
and replacing it by one of the other trivalent trees from 
Figure~\ref{pic_mutation}.
Two graphs are \textbf{mutation-equivalent} 
if they may be transformed by sequence of mutations into each other.
\end{defin}
\begin{rmk}
By definition it is not possible to mutate 
along an edge that forms a cycle of length one.
However, for longer cycles, mutations are possible, 
and one shortens the length of the cycle.
Figure~\ref{pic_mutation_along_cycle_edge} shows this fenomena
for Hammock graph mutation equivalent to LittleMan.

\myVSPACEfigure\begin{center}
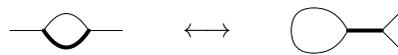
  \begin{minipage}{\textwidth}
\[
\begin{xy}
(0,0);(0,1.5):;
(0,-5)*{} ; (0,-2)*{} **@{-};
(0,-2)*{} ; (0, 2)*{}  **[|<1.5pt>]\crv{~**\dir{-} (-3,0)};
(0,-2)*{} ; (0, 2)*{}  **\crv{~**\dir{-} ( 3,0)};
(0, 5)*{} ; (0, 2)*{} **@{-};
\end{xy} 
\qquad \longleftrightarrow \qquad
\begin{xy}
(0,0);(0,1.5):;
( 0, 0)*{} ; ( 2, 3)*{} **\crv{~**\dir{-} (2,1)};
( 2, 3)*{} ; ( 0, 5)*{} **\crv{~**\dir{-} (2,5)};
(-2, 3)*{} ; ( 0, 5)*{} **\crv{~**\dir{-} (-2,5)};
( 0, 0)*{} ; (-2, 3)*{} **\crv{~**\dir{-} (-2,1)};
( 0, 0)*{} ; ( 0,-3)    **[|<1.5pt>]@{-};
( 0,-3)*{} ; (-2,-5)*{} **@{-};
( 0,-3)*{} ; ( 2,-5)*{} **@{-};
\end{xy}
\]
\captionof{figure}{Mutation along cycle edge shortens a cycle}
\label{pic_mutation_along_cycle_edge} 
 \end{minipage} \end{center}\myVSPACEfigure
\end{rmk}

\begin{lem}\label{cycle_shorter}
Suppose edges $\{e_1,\ldots,e_k\}$ form a cycle in the graph $\ccG$
and assume $k>1$. 
Then, for any $i\in\{1,\ldots,k\}$ one of the two mutations 
along $e_i$ shortens the cycle $\{e_1,\ldots,e_k\}$ by one 
in the resulting graph, i.e.~$\{e_1,\ldots,e_{i-1},e_{i+1},\ldots,e_k\}$ 
is a cycle in the new graph.
\end{lem}

A \textbf{caterpillar} is a trivalent tree, which after removing 
all leaves and petioles becomes a string of edges
as shown on Figure~\ref{pic_caterpillar_tree}.

\myVSPACEfigure\begin{center}
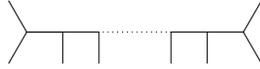

\begin{minipage}{\textwidth}
\[
\begin{xy}
=(  8, 0     )"A" ;
=(  4, 6.9282)"B" ;
=( -4, 6.9282)"C" ;
=( -8, 0     )"D" ;
=( -4,-6.9282)"E" ;
=(  4,-6.9282)"F" ;
=( 16, 0)"Aleg";
=(-16, 0     )"Dleg";
=(-20, 6.9282)"Hone";
=(-20,-6.9282)"Htwo";
<0mm,0mm>;<0.6mm,0mm>:;
"Hone"; "Dleg"       **@{-};
"Htwo"; "Dleg"       **@{-};
"Dleg" ; "D"         **@{-};
"Dleg"; (0,0)        **@{-};
(0, 0); (16,0)       **@{.};
(16,0); (32,0)       **@{-};
"D"   ; (-8,-6.9282) **@{-};
(0, 0); ( 0,-6.9282) **@{-};
(16,0); (16,-6.9282) **@{-};
(24,0); (24,-6.9282) **@{-};
(32,0); (36,-6.9282) **@{-};
(32,0); (36, 6.9282) **@{-};
\end{xy}
\]  
\captionof{figure}{Caterpillar tree}
\label{pic_caterpillar_tree}
\end{minipage}
\end{center}\myVSPACEfigure

\begin{lem}\label{graphs_mutation_equivalent}
Let $\ccG_1$ and $\ccG_2$ be connected, trivalent graphs 
both with $n$ leaves and first Betti numbers $g$.
Then they are mutation-equivalent.
Moreover, for any ordered subsets $S_1$ of cycle edges of $\ccG_1$
and $S_2$ of cycle edges of $\ccG_2$, of the same size,
both with the property that 
removing $S_i$ from $\ccG_i$ does not disconnect the graph,
we can find a sequence of mutations that 
avoid the edges from $S_1$ ($S_2$) and
sends $i^{th}$ edge of $S_1$ to  the $i^{th}$ edge of $S_2$.
Also, any mutation sends a leaf of $\ccG_1$ to a leaf of $\ccG_2$.
\end{lem}
\begin{prf}
Let $\ccG$ be connected, trivalent graph with $n$ leaves and 
the first Betti number $g$ and $S$ a subset of cycle edges as above.
We will prove that $\ccG$ is mutation-equivalent to 
a trivalent graph obtained by attaching $g$ cycles of length one
to a caterpillar tree with $n+g$ leaves.
We will choose  mutations so that they will satisfy the required property.

\emph{Step 1.} 
We proceed by induction on the set $S$ and the first Betti number of $\ccG$.
For an edge $e \in S$ we can find a cycle that contains 
no other elements of $S$. This is because after removing
all edges from $S$ the graph $\ccG$ is connected, so there
is a  path $\gamma$ from $\partial_1(e)$ to $\partial_2(e)$, which 
together with $e$ form the required cycle.
By repeatedly using Lemma~\ref{cycle_shorter},
we reduce the length of this cycle to one, by
performing mutations along edges from  $\gamma$.
In the new graph the edge $e$ forms a loop.
We can consider this graph with $e$ removed,
reducing both the size of $S$ and $g$.
If $|S|<g$, then we repeat the above $g-|S|$ times
starting from any cycle edge, which is not a loop.

After repeating this procedure $g$ times,
we get  a tree with $g$ loops
(all  edges from $S$ are among them) attached to some leaves.
We can assume that this tree is a caterpillar,  
as we know from \cite[Lem.~2.18]{buczynska_wisniewski}, 
that any trivalent tree is mutation-equivalent to a caterpillar 
with the same number of leaves.

\emph{Step 2.}
We observe that it does not matter to which leaves the cycles are attached,
we can move a cycle from a leaf to any another leaf.

In Figure~\ref{pic_mut_equivalent} we illustrate those two steps.

\myVSPACEfigure\begin{center}  \begin{minipage}{\textwidth}
\[
\begin{xy}
=(  8, 0     )"A" ;
=(  4, 6.9282)"B" ;
=( -4, 6.9282)"C" ;
=( -8, 0     )"D" ;
=( -4,-6.9282)"E" ;
=(  4,-6.9282)"F" ;
=( 16, 0)"Aleg";
=(-16, 0     )"Dleg";
=(-20, 6.9282)"Hone";
=(-20,-6.9282)"Htwo";
=(-19.6, 7.2282)"HoneXbase";
=(-20.4,-6.6282)"HtwoXbase";
<0mm,0mm>;<0.6mm,0mm>:;
\save
(0,12)*{}; (0,-12)*{};(0,0); 
\restore
\drop{\LEGS------}
\drop{\HEX------}
\POS{"Dleg"+"C"} ; "Dleg" **@{-};
\POS{"Dleg"+"E"} ; "Dleg" **@{-};
\save
"Hone";"HoneXbase":0*{\CHERRY};
\restore \save
"Aleg" ; (26, 0)  **\crv{~**\dir{-} (21, 10)};
"Aleg" ; (26, 0)  **\crv{~**\dir{-} (21,-10)};
(26, 0) ; (34,0) **@{-};
\restore \save
{\ar^{\scriptstyle{Step\ 1}}@{~>} (40,0) ; (60,0)*{}};
\restore \save
(85,0);(86,0):;
\save
\POS{"Hone"};\POS{"HoneXbase"}:0*{\CHERRY};
\restore \save
"Hone"; "Dleg" **@{-};
"Htwo"; "Dleg" **@{-};
"Dleg" ; "D" **@{-};
"Dleg"; (32,0) **@{-};
"D" ;   (-8,-6.9282) **@{-};
(0, 0); ( 0,-6.9282) **@{-};
\restore \save
( 0,-6.9282) ; (-0.5, -6.9282):0*{\CHERRY};
\restore \save
(8, 0); ( 8,-6.9282) **@{-};
(16,0); (16,-6.9282) **@{-};
(24,0); (24,-6.9282) **@{-};
\restore \save
( 24,-6.9282) ; (23.5, -6.9282):0*{\CHERRY};
\restore \save
(32,0); (36,-6.9282) **@{-};
(32,0); (36, 6.9282) **@{-};
\restore \save
\restore \save
{\ar^{\scriptstyle{Step\ 2}}@{~>} (40,0) ; (60,0)*{}};
(85,0);(86,0):;
\save
\POS{"Hone"};\POS{"HoneXbase"}:0*{\CHERRY};
\restore \save
\POS{"Htwo"};\POS{"HtwoXbase"}:0*{\CHERRY};
\restore \save
"Hone"; "Dleg" **@{-};
"Htwo"; "Dleg" **@{-};
"Dleg" ; "D" **@{-};
"Dleg"; (32,0) **@{-};
"D" ;   (-8,-6.9282) **@{-};
\restore \save
( -8,-6.9282) ; (-8.5, -6.9282):0*{\CHERRY};
\restore \save
(0, 0); ( 0,-6.9282) **@{-};
(8, 0); ( 8,-6.9282) **@{-};
(16,0); (16,-6.9282) **@{-};
(24,0); (24,-6.9282) **@{-};
(32,0); (36,-6.9282) **@{-};
(32,0); (36, 6.9282) **@{-};
\restore 
\end{xy}
\]   

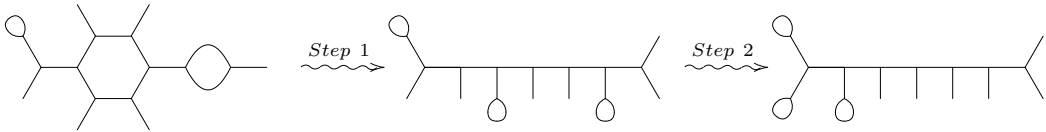
\captionof{figure}{Every graph is mutation-equivalent to a caterpillar graph 
\label{pic_mut_equivalent}
}
 \end{minipage} \end{center}\myVSPACEfigure
The last claim  follows simply form the definition:
 mutation maps an inner edge to an inner edge, and a leaf to a leaf.
\end{prf}

\subsection{Definition as quotient.}\label{section_def_of_model}

Given a not necessarily connected trivalent graph $\ccG$, 
we construct a toric variety $X(\ccG)$, 
generalizing the binary symmetric model of trivalent tree 
from \cite{buczynska_wisniewski}.

As we have already explained, see equation~\eqref{graph_from_tripods}, 
any trivalent graph is the union of $|\ccV|-n$
elementary tripods with some edges identified.
To define the variety $X(\ccG)$ we replace each elementary tripod $\TRIPOD_v$ 
with $\PP^3_v$, union with product, 
and the edge identification with
a quotient by an action of a one-parameter torus.

\begin{defin}\label{defin_as_quotient}
Let $\ccG$ be a trivalent graph.
To an inner vertex $v \in \ccV \setminus \ccL$ 
we associate projective space $\PP^3_v$
with coordinates $x^v_{\emptyset},x^v_{12},x^v_{13},x^v_{23}$.
To any edge $e\in\ccE$ we associate an action  $\lambda^e_v$ of $\CC^*$ on $\PP^3_v$ 
with weights $0$ and $1$ as follows:
\[
\lambda^e_v(t)(x_S)=
\begin{cases}
  t \cdot x^v_S & \text{if the index of $i_v^{-1}(e)\in\{e_1,e_2,e_3\}$ belongs to the set $S$,}\\
  x^v_S         & \text{otherwise.}
\end{cases}
\]
So we have an action of a three-dimensional torus on 
$\PP^3=\Proj \CC [x^v_{\emptyset},x^v_{12},x^v_{13},x^v_{23}]$ 
with weights:
\[
\begin{bmatrix}
0 & 1 & 1 & 0 \\
0 & 1 & 0 & 1 \\
0 & 0 & 1 & 1 \\
\end{bmatrix}
\]
In other words if for example $i_v(e_2)=e$ then 
$\lambda^e_v$ acts with weight $1$ on $x^v_{12}$ and $x^v_{23}$
and with weight $0$ on $x^v_{\emptyset}$ and $x^v_{13}$.
This action extends to an action $\lambda^e_v$ on 
$\prod_{v \in \ccV \setminus \ccL} \PP^3_v$ which is non-trivial 
only if $v$ is an end of the edge $e$.
Thus, for any internal edge $e\in \ccE\setminus \ccP$, we can define
a $\CC^*$-action $\lambda^e_{\partial_1(e)} \times -\lambda^e_{\partial_2(e)}$ on 
\begin{equation}\label{product_of_P3s}
\prod_{v \in \ccV \setminus \ccL} \PP^3_v   
\end{equation}
to be the product action of the action $\lambda^e_{\partial_1(e)}$ on $\PP^3_{\partial_1(e)}$ and
the action $\lambda^e_{\partial_2(e)}$ with  opposite weights on $\PP^3_{\partial_2(e)}$.
We define the \textbf{ phylogenetic model of a trivalent graph} $\ccG$ 
to be the good quotient:
\begin{equation}\label{model_as_quotient}
X(\ccG) := 
\Bigl(
\prod_{v \in \ccV \setminus \ccL} \PP^3_v 
\Bigr)
\quad \biggit  \prod_{e\in\ccE\setminus\ccP} 
\left({\lambda^e_{\partial_1(e)} \times -\lambda^e_{\partial_2(e)}} \right )
\end{equation}
of the toric variety $ \prod_{v \in \ccV \setminus \ccL} \PP^3_v $
by a subtorus of dimension $|\ccE|-|\ccL|$ of the torus.
The subtorus by which we are dividing is a product of all the $\CC^*$'s
over all internal edges of the graph $\ccG$ and 
the linearized line bundle is $\boxtimes_{v \in \ccV(\ccG)} \ccO_{\PP^3_v}(1)$.
By Theorem~\ref{toricGIT}.
$X(\ccG)$ is toric as it is the quotient of a toric variety by a subtorus.
\end{defin}

\begin{rmk}
In Definition~(\ref{model_as_quotient}) the choice that we made
defining the action of the torus $(\CC^*)^{\ccV \setminus \ccL}$ 
only depends on the choice of coordinates of the torus.  
If we choose different orientation of the edge $e$, then
the two $\CC^*$-actions $\lambda^e_{\partial_1(e)} \times -\lambda^e_{\partial_2(e)}$ 
and $-\lambda^e_{\partial_1(e)} \times \lambda^e_{\partial_2(e)}$
differ by composing with $t \mapsto {\frac 1 t}$.
\end{rmk}

\begin{rmk}\label{torus_of_leaves}
Let $l$ be a leaf of a graph $\ccG$ and $e$ the adjacent petiole.
The action $\lambda^e_l$ descends to a non-trivial action on 
the quotient variety $X(\ccG)$ and is denoted by $\lambda^l$.
For a subset  $\ccS\subset\ccL$ of the leaves of $\ccG$ 
of cardinality $k$ we have an action of a $k$-dimensional torus 
$\TT(\ccS)$ --- a product of the corresponding $\lambda^l$'s.
\end{rmk}

If we set $k=|\ccV|-n$ the number of inner vertices,
we can rewrite $|\ccV|-|\ccE|=|\comp\ccG|-g$ using 
Lemma~\ref{graph_invariants}\ref{trivalent_equation} 
to get $k = (2g-2|\comp\ccG|) + n$. 
Observe that
\begin{itemize}
\item  $k$ is the number of the $\PP^3$'s in~\eqref{product_of_P3s}
\item $g$ is the first Betti number.
\end{itemize}
On the other hand, $|\ccE| - n = k + g - |\comp\ccG|$
is the number of inner edges which is the dimension of the torus 
that we divide by in~\eqref{model_as_quotient}.
We get a variety of dimension 
\[
\dim X(\ccG) = 3k - (|\ccE|-n) = 3g - 3|\comp\ccG| + 2n = |\ccE|.
\]


\subsection{Lattice and cone.}\label{section_def_lattice_and_cone}
Given a trivalent graph $\ccG$ we construct the toric data that
allows to recover its toric model. The graded lattice, denoted by $\Mg$,
as well as the graded cone in it
have both rank one bigger than the dimension of the model $X(\ccG)$, 
and the latter is equal to the number of edges $\ccE(\ccG)$.

\begin{defin}\label{lattice_defin}
Given a graph $\ccG$ let $\ZZ \ccE=\bigoplus \ZZ \cdot e$ 
be the lattice spanned by $\ccE$,
and $\ZZ \ccE^{\vee} = \Hom (\ZZ \ccE, \ZZ)$ be its dual.
Elements of the lattice $\ZZ \ccE$ are formal linear combinations of the edges, 
forming the standard basis of~$\ZZ \ccE$. The dual lattice~$\ZZ \ccE^{\vee}$ 
comes with the dual basis $\{ e^* \}_{e\in \ccE}$.
We identify vertices of the graph $\ccG$ 
with certain elements of~$\ZZ \ccE^{\vee}$: 
\begin{equation}\label{vertex_rep_in_N}
v=\sum_{e \ni v}e^*.
\end{equation}
We also define 
$M = \{u\in \ZZ \ccE :\forall v\in\ccV\ v(u)\in 2\ZZ\}$
and its dual $N=\Hom(M,\ZZ)$.
Then the graded lattice of the graph is
\[
\Mg=\ZZ \oplus M, 
\]
with the degree map
\[\deg : \Mg = \ZZ \oplus M \ra \ZZ,\]
which is the projection to the first summand.
The \textbf{degree} of $\omega \in \Mg$ is $\deg(w)$.
\end{defin}
If there is more then one graph in question we will write 
$M(\ccG)$ for $\Mg(\ccG)$ etc.

Let us use the following notation for the elements of the lattice 
$\ZZ \ccE^{\vee}$ dual to the edges meeting at the vertex $v$
\[
\begin{array}{lcr}
a_v:=\big(i_v(e_1)\big)^*,  & 
b_v:=\big(i_v(e_2)\big)^*,  & 
c_v:=\big(i_v(e_3)\big)^*,  \\
\end{array}
\]
where $\{e_1, e_2, e_3\}$ are the edges of $\TRIPOD$ and 
$i_v : \TRIPOD \mono \ccG$ is, as before, a map which is
locally an embedding and sends the central vertex of the $\TRIPOD$
to $v$ --- an inner vertex of~$\ccG$.
Whenever we use this notation we have a fixed presentation
as in~(\ref{graph_from_tripods}).

Given an element $\omega$ in either $\ZZ \ccE$, $M$ or $\Mg$,
each of $a_v, b_v, c_v \in \ZZ \ccE^{\vee}$ measures the coefficient of $\omega$
at an edge incident to $v$.
Then~\eqref{vertex_rep_in_N} becomes
\[
v=a_v+b_v+c_v.
\]

\begin{defin}
The \textbf{degree} of $\omega\in \Mg$ \textbf{at a vertex} $v\in\ccV(\ccG)$ is
\[
\deg_v(\omega):=\half\cdot\bigl(a_v(\omega)+b_v(\omega)+c_v(\omega)\bigr).
\]
The \textbf{minimal degree} of $\omega$ is
\[\deg_{\min}(\omega):= \max_{v\in\ccV}\{\deg_v(\omega) \}, \]
where $ \pi_M : \Mg \ra M$ is the projection to the second summand.
\end{defin}
The name minimal degree will be clear after we define the cone $\tau(\ccG)$.

We identify paths and networks in $\ccG$ with 
elements of the lattices $M$ and $\Mg$, by replacing 
union with sum in the group $\ZZ \ccE$
\begin{defin} 
A \textbf{network in the graded lattice} $\Mg$ is a pair
$\omega=(1,a) \in \Mg$ where $a\in M$ is a network.
\end{defin}

\begin{lem}
An element of the lattice $M$ is represented by a labeling 
of the edges of $\ccG$ with integers so that 
the sum at any vertex is even. 
Thus the lattice $M \subset \ZZ \ccE$ is generated by
\begin{enumerate}
\item networks
\item $\{ 2e \mid e \in \ccE \}$.
\end{enumerate}
\end{lem}
\begin{prf}
Let $\omega \in M$. By using generators of the second type, 
we can assume that  $0 \leq a_v(\omega),b_v(\omega),c_v(\omega) \leq 1$ for any vertex $v$. 
This implies that $a_v(\omega)+b_v(\omega)+c_v(\omega)$ equals $0$ or $2$ 
and either two among $a_v(\omega),b_v(\omega),c_v(\omega)$ are one or all are zero.
This implies that $\omega$ is a network, since it corresponds to 
a disjoint union of path. A path goes through a vertex $v$
means in terms of $a_v(\omega),b_v(\omega),c_v(\omega)$ that exactly two of them 
are one.
\end{prf}

We define the cone $\tau(\ccG)$ of the graph, which is 
the semigroup defining the model of the graph as
projective spectrum of the semigroup ring, 
as we will see in Theorem~\ref{model_as_proj}.

\begin{defin}\label{def_cone_of_G}
For a graph $\ccG$ we define its cone  $\tau = \tau(\ccG) \subset  \Mg$ 
as the set of $\omega \in \Mg$ which satisfy following inequalities:
\begin{enumerate}
\item $a_v(\omega),b_v(\omega),c_v(\omega) \ge 0$,
\item for any vertex $v\in \ccV$ triangle inequalities hold
\[
|a_v(\omega)-b_v(\omega)| \leq c_v(\omega) \leq a_v(\omega)+ b_v(\omega),\text{ and}
\]
\item \label{mindeg}$\deg(\omega) \geq \deg_{\min}(\omega)$.
\end{enumerate}
\end{defin}

\begin{rmk}
To explain the name minimal degree $\deg_{\min}(\omega)$, note that 
for any $\omega$ in the cone we have the following equality
\[
\deg_{min}(\omega)=\min_{\omega' \in \tau} \{\deg(\omega') : \pi_{\hM}(\omega')=\pi_{\hM}(\omega)\}.
\]
\end{rmk}
\begin{prf}[Proof of the remark]
By part~\ref{mindeg} of the Definition~\ref{def_cone_of_G} of $\tau$ 
for any $\omega' \in \tau$ satisfying $\pi_{\hM}(\omega')=\pi_{\hM}(\omega)$  we have 
\[
\deg(\omega') \ge \deg_{\min}(\omega') = \deg_{\min}(\omega) 
\] 
since by definition of $\deg_{\min}(\omega)$ only depends on $\pi_{\hM}(\omega)$.
This means 
\[
\deg_{\min}(\omega) \le \min_{\omega' \in \tau} \{\deg(\omega') : \pi_{\hM}(\omega')=\pi_{\hM}(\omega)\}.
\]
To prove the equality we will find $\omega' \in \tau$ 
with $\deg(\omega')=\deg_{\min}(\omega)$.
Let us write $\omega=(\deg(\omega), \alpha)$ as it is an element of $\Mg$.
Recall that  $\deg_{\min}(\omega)$ is the maximum of 
$\deg_v(\omega)=\half \cdot (a_v+b_v+c_v)(w)$ over all vertices $v$ of the graph. 
Thus $\omega'=(\deg_{\min}(\omega),\alpha)$ is in the cone $\tau$ and 
has the required degree in $\Mg$.
\end{prf}

\begin{thm}\label{model_as_proj}
The variety $X(\ccG)$ is isomorphic to the toric variety $\Proj \CC[\tau(\ccG)]$.
\end{thm}
\begin{prf}
To see this we first observe that each $\PP^3_v$ in 
Definition~\ref{defin_as_quotient} of the model $X(\ccG)$ can be written as 
$\Proj \CC [\tau(\TRIPOD_v)]$, where each cone
\[
\tau(\TRIPOD_v)=\conv(0000,1000,1110,1101,1011)
\]
is a cone over a tetrahedron and is clearly defined by 
the required inequalities.
Next, taking the product of $\PP_v$ corresponds by 
Lemma~\ref{product_of_varieties_is_a_graded_product_of_cones}
to taking a graded product of cones. 
Thus the product cone is defined by required inequalities.
Lastly we use the description of the quotient of a toric variety
by a subtorus of the torus in Theorem~\ref{toricGIT}.
The subtorus in Definition~\ref{defin_as_quotient}
is a product of the $\CC^*_e$ over all internal edges.
Taking the quotient with respect to such a torus
corresponds to cutting the cone with the hyperplane 
of the type $a_{\partial_1(e)}=b_{\partial_2(e)}$, which preserves the inequalities.
\end{prf}

\subsection{$\ZZ$-generators of the cone $\tau(\ccG)$}
\label{section_cone_gens}
Knowing that the model $X(\ccG)$ is the projective spectrum
of the semigroup algebra of $\tau(\ccG)$
means that it is a subvariety of a weighted projective space
with weights equal to the degrees of the chosen generators.
When $\ccG = \ccT$ is a tree the cone is generated in degree~$1$
so the embedding is into a (straight) projective space
$\PP^k=\PP(1, \ldots, 1)$. 
In this case, by \cite{buczynska_wisniewski}, we already know 
all about this cone, see Proposition~\ref{tree_gens} below.
It is represented by its degree $1$ section --- a normal lattice polytope
$\Delta(\ccT) \subset {1} \times M$, whose vertices span
the cone $\tau(\ccT)$ and the lattice points generate the semigroup. 
Our goal is to show that for graphs with the first Betti number one,
the semigroup is generated in degrees~$1$ and~$2$.

In order to describe $\ZZ$-generators of the cone $\tau(\ccG)$
we will  express elements of $\tau(\ccG)$ in terms of $\ccG$. 
We also decompose the graph $\ccG$ into smaller graphs
for which $\ZZ$-generators of the corresponding cones 
are easier to find. 

We explain that any element $\omega$ of the cone $\tau(\ccG)$ 
locally decomposes into paths. 
In the graph $\TRIPOD$ there are three non-empty paths,
each consisting of two edges. Let us denote them by:
\[\begin{array}{lcr}
x:=e_2+e_3, & y:=e_1+e_3,&z:=e_1+e_2,\\ 
\end{array}\]
where $\{e_1,e_2,e_3\}$ are edges of $\TRIPOD$.

For an arbitrary $\ccG$, we know that
given an element $\omega \in \tau(\ccG)$ of the cone  and a vertex $v \in \ccV(\ccG)$, 
the numbers $a_v(\omega), b_v(\omega), c_v(\omega)$ satisfy the triangle inequalities
and their sum is even.
This allows us to, locally at $v$, rewrite $\omega$ 
as sum of paths $x$, $y$ and~$z$. 
The picture of this decomposition is drawn in Figure~\ref{local_paths}.
\myVSPACEfigure\begin{center}
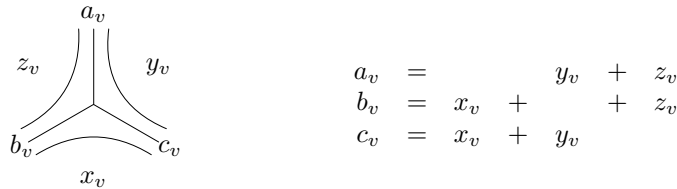
  \begin{minipage}{\textwidth}
\[
\begin{xy}
(0,0)*{}="0";
(0,10)*{} ="A"; 
(-2,0)="Abvector";
(1,-\rootthree)="Bcvector";
(1, \rootthree)="Cavector";
( 1,10)*{} ="Ac";
(-\fiverootthree,-5)*{}="B" ; 
( \fiverootthree,-5)*{}="C" ; 
(0,-2)="a"; ( \halfrootthree,0.5)="b"; (-1.5,  \halfrootthree)="c";
"0" ; "A" **@{-}-"a"*{a_v};
"0" ; "B" **@{-}-"b"*{b_v};
"0" ; "C" **@{-}-"c"*{c_v};
"A"+"Abvector" ; "B"-"Bcvector" **\crv{"0"+"c"};
"B"+"Bcvector" ; "C"-"Cavector" **\crv{"0"+"a"};
"C"+"Cavector" ; "A"-"Abvector" **\crv{"0"+"b"};
"0"-"A" *{x_v};
"0"-"B" *{y_v};
"0"-"C" *{z_v};
\end{xy}
\qquad\qquad\qquad
\begin{array}{rlrlrlrl}
a_v&=&     & & y_v &+& z_v\\
b_v&=& x_v &+&     &+& z_v\\
c_v&=& x_v &+& y_v & &    \\
\end{array}
\]
\captionof{figure}{Local paths around a vertex}
\label{local_paths}
 \end{minipage} \end{center}\myVSPACEfigure

Our aim is to find the $\ZZ$-generators by understanding how the graph $\ccG$
was built from smaller pieces. 
Each of the operations in Definition~\ref{definition_operations_on_graphs}  
has a corresponding operation on lattices and cones. 
By Definition~\ref{defin_as_quotient}, the model of a disjoint sum 
of graphs is the product of the models,
so the underlying cone is the graded product of corresponding cones.

\begin{lem}\label{cone_of_graft_is_fiberproduct}
Let $\ccG_1$ and $\ccG_2$ be two trivalent graphs then  
\[
\begin{array}{lcl}
\Mg(\ccG_1 \sqcup \ccG_2)&=&\Mg(\ccG_1)\times_{g}\Mg(\ccG_2)\text{, and} \\
 \tau(\ccG_1 \sqcup \ccG_2)&=& \tau(\ccG_1)\times_{g} \tau(\ccG_2). \\
\end{array}
\]
\end{lem}

In the definition of $X(\ccG)$ we take a quotient by a torus
corresponding to the set of inner edges.
In other words we have translated the operation of gluing of two leaves of a graph $\ccG$ 
into taking a quotient by appropriate $\CC^*$-action of  the model $X(\ccG)$.
The following observation says that if we choose to glue 
some pairs of leaves first and then the rest of the pairs
it does not matter how we partition the set of pairs of leaves 
or which order we choose. In all cases the resulting variety is the same.

\begin{lem}\label{quotient_commutative}
Let two reductive, commutative groups $H_1$ and $H_2$ and their direct sum 
\mbox{$H_1 \oplus H_2$} act  on a projective variety $X$. 
Suppose all those actions are linearized with respect to 
some ample Weil divisor $\ccL$. Then 
\[
X \git (H_1\oplus H_2)=(X \git H_1) \git H_2=(X \git H_2) \git H_1,
\]
where the semistable points on $X$ are taken with respect to $\ccL$ and 
on quotients of $X$ with respect to the push-forward of $\ccL$.
\end{lem}

We have identified the vertices of $\ccG$ with elements of $\ZZ \ccE^{\vee}$, 
see \eqref{vertex_rep_in_N}. We observed in Remark~\ref{torus_of_leaves} 
that a leaf $l$  yields a $\CC^*$-action $\lambda^l$ on $X(\ccG)$. 
Given two leaves $l_1$ and $l_2$ of $\ccG$, by definition we have
\[
X(\ccG^{l_1}_{l_2} \glue) = X(\ccG) \git \left(\lambda^{l_1} 
\times -\lambda^{l_2}\right).
\]
In terms of toric geometry this quotient corresponds to 
the intersection of $\tau(\ccG)$ with the kernel of $l_1-l_2$, 
where we  treat $l_1$ and $l_2$ as elements of the lattice $(\Mg)^{\vee}$.
Thus the following lemma is a consequence of Theorem~\ref{toricGIT}.
\begin{lem}
Let  $l_1$ and $l_2$ be distinct leaves of $\ccG$. Then
\[
\begin{array}{lclcl}
\Mg(\ccG^{l_1}_{l_2}\glue)&=&\Mg(\ccG)&\cap& \ker (l_1 - l_2)\\
 \tau(\ccG^{l_1}_{l_2}\glue )&=& \tau(\ccG)&\cap& \ker (l_1 - l_2).
\end{array}
\]
\end{lem}

The reverse operation on graphs is to cut an edge into two new edges.
In the next lemma we explain how this is reflected on the cones and lattices.

\begin{lem}\label{cut_edge}
Let $\ccG$ be a trivalent graph and $e \in \ccE \setminus \ccP$ 
an internal edge. As before $\ccG^e$ is the graph obtained from $\ccG$
by cutting the edge $e$.
Then there are natural embeddings of the cones and graded lattices:
\begin{align*}
i^e &: \Mg(\ccG) \mono \Mg(\ccG^e) & i^e &: \tau(\ccG) \mono \tau(\ccG^e)  
\end{align*}
\end{lem}
\begin{prf}
Let 
\[
\ZZ \ccE^e(\ccG)= \bigoplus_{e' \in \ccE(\ccG)\setminus \{e\}} \ZZ \cdot e'
\]
be the lattice spanned by all other edges.
We can decompose the lattices $(\ZZ \ccE)^{gr}(\ccG)$ 
and $(\ZZ \ccE)^{gr}(\ccG^{e})$ 
\begin{align*}
\ZZ \ccE(\ccG)   &=\ZZ\oplus \ZZ \ccE^e(\ccG) \oplus \ZZ \cdot e,\\
\ZZ \ccE(\ccG^e) &=\ZZ\oplus \ZZ \ccE^e(\ccG) \oplus \ZZ \cdot  e_1 \oplus \ZZ\cdot e_2.
\end{align*}
Thus we can embed the lattices identifying the first two summands and 
taking a diagonal embedding of the third one:
$i^e :=\id_{\ZZ}\oplus \id_{\ZZ \ccE^e} \oplus \Delta_{\ZZ} : 
(\ZZ \ccE)^{gr}(\ccG) \mono (\ZZ \ccE)^{gr}(\ccG^{e})$.
To complete the proof we check that $i^e$ restricts to 
the lattices with parity condition
and to the cones.
\begin{align*}
\xymatrix{
  \ZZ \ccE(\ccG)  \ar@{^(->}[r]^{i^e}               & \ZZ \ccE(\ccG^e) \\
  \Mg(\ccG) \ar@{^(->}[r]^{i^e} \ar@{^(->}[u] & \Mg \ar@{^(->}[u](\ccG^e)\\
}&&
\xymatrix{
  \ZZ \ccE(\ccG)  \ar@{^(->}[r]^{i^e}               & \ZZ \ccE(\ccG^e) \\
  \tau(\ccG) \ar@{^(->}[r]^{i^e} \ar@{^(->}[u] & \tau \ar@{^(->}[u](\ccG^e)\\
}&
\end{align*}
To see it we only need to check the parity condition about the ends 
$\partial_1(e)$ and $\partial_2(e)$ of the edge $e$ that we cut. 
By definition $e_1^*(i^e(e)) = e_2^*(i^e(e))$, so 
for any $\omega \in \Mg$ we have 
$\partial_1(e)(\omega) = \partial_1(e_1)(i^e(\omega))$ and $\partial_2(e)(\omega) = \partial_2(e_2)(i^e(\omega))$.
In the same way $i^e$ preserves all the inequalities defining 
cones $\tau(\ccG)$ and $\tau(\ccG^e)$.
\end{prf}

When the edge $e$ is not a cycle edge, the graph $\ccG^e$ is not connected.
We write $\ccG^e = \ccG_1 \sqcup \ccG_2$, where $\ccG_1$ (respectively $\ccG_2$) 
is the part containing $e_1$ (respectively $e_2$). 
Then we have a projection $\pi^e_1$ (respectively $\pi^e_2$) of lattices
\[
\pi^e_1 : \Mg(\ccG^e)=\Mg(\ccG_1) \times \Mg(\ccG_2) \epi \Mg(\ccG_1)
\]
which restricts to a projection of cones.
For a non-cycle edge $e \in \ccE$   we denote by $\rho_1^e$ 
(respectively $\rho_2^e$) the composition $\rho_1^e= \pi_1^e \circ i^e $  
of the above defined maps.

\begin{rmk}
When $e\in\ccE$ is not a cycle edge we write $\ccG^e=\ccG_1 \sqcup \ccG_2$.
Then the cone $\tau(\ccG)$ is the following fiber product 
of the cones $\rho_1^e(\tau(\ccG))=\tau(\ccG_1)$ and  
$\rho_2^e(\tau(\ccG))=\tau(\ccG_2)$.
The same is true for the lattice $\Mg(\ccG)$.
\begin{align*}
\xymatrix{
\tau(\ccG)   \ar[r]^{\rho_1^e} \ar[d]^{\rho_2^e}& 
\tau(\ccG_1) \ar[d]^{\deg \oplus e_1^*}\\
\tau(\ccG_2) \ar[r]^{\deg \oplus e_2^*} & \ZZ \oplus \ZZ
}&&
\xymatrix{
\Mg(\ccG)   \ar[r]^{\rho_1^e} \ar[d]^{\rho_2^e}& \Mg(\ccG_1) 
\ar[d]^{\deg \oplus e_1^*}\\
\Mg(\ccG_2) \ar[r]^{\deg \oplus e_2^*} & \ZZ \oplus \ZZ
}&
\end{align*}
\end{rmk}
Now we turn our attention to our main task of finding the $\ZZ$-generators 
of the cone $\tau(\ccG)$.
\begin{lem}
For any graph $\ccG$ the set of degree $1$ integer points of cone $\tau(\ccG)$ 
is equal to the set of networks.
\end{lem}
\begin{prf}
If $\omega \in \tau$ is a point in the cone of degree $1$, then for any vertex $v \in \ccV$,
\[
1=\deg(\omega) \ge \deg_{\min}(\omega) = \max_{u\in \ccV} \{ \deg_u (\omega) \} \geq \deg_v (\omega) \ge 0.
\]
By definition $\deg_v= x_v+y_v+z_v \ge 0$ so exactly one of
$x_v(\omega)$, $y_v(\omega)$, $z_v(\omega)$ equals one and the other two are zero
or all are zero.
Equivalently exactly two of $a_v(\omega)$, $b_v(\omega)$, $c_v(\omega)$ are one,
and the third one is zero, or all are zero. 
This means that $\omega$ is a network.
\end{prf}
\begin{cor}  
All networks are among the minimal $\ZZ$-generators of the cone $\tau(\ccG)$.
\end{cor}
In fact when the graph in question is a  tree these are the $\ZZ$-generators.
\begin{prop}[{\cite[\S 2.1]{buczynska_wisniewski}}]
\label{tree_gens}
If $\ccG$ is a trivalent tree, 
then $\tau(\ccG)$ is generated in degree~$1$. 
Moreover, the generators are exactly networks of paths, 
which in this case are determined by their values on the leaves.
Thus a generator of $\tau(\ccG)$ is identified with 
a sequence of $0$'s and $1$'s of length $n(\ccG)$ 
with even number of $1$'s.
As a consequence a model of a trivalent tree with $n$ leaves 
comes with an embedding into projective space~$\PP^{2^{n-1}-1}$.
\end{prop}

Let $\ccG$ be a graph with the first Betti number one.
We will describe the generators of the semigroup $\tau(\ccG)$ 
in this case. We cut all the cycle legs $l_1,\ldots,\l_k$ of $\ccG$ and write
$\ccG^{l_1,\ldots,l_k}=\ccG_0 \sqcup\ccG_1 \sqcup\ldots \sqcup\ccG_k$, where
$\ccG_0$ is a polygon graph and thus $\ccG_1,\ldots, \ccG_k$ are pendant trees 
(see Definition~\ref{definition_cycle_edge_pending_tree_polygon_graph}).
Thus any element $\omega \in \tau(\ccG)$ has a lift 
$\tilde\omega \in \tau(\ccG_0 \sqcup \ldots \sqcup \ccG_k)
= \tau(\ccG_0) \times_g \tau(\ccG_1) \times_g \ldots \times_g \tau(\ccG_k)$
and \textbf{components} $\tilde\omega=(\omega_0,w_1,\ldots,\omega_k)$, 
which can be written $\omega_i = \rho^{l_i}(\omega)$.

\begin{thm}\label{cone_gens}
Let $\ccG$ be a trivalent graph with the first Betti number exactly one.
Any element $\omega \in \Mg(\ccG)$ is a minimal $\ZZ$-generator of $\tau(\ccG)$
if and only if it satisfies one of the following conditions
\begin{enumerate}
\item $\omega$ has degree $1$ and $\omega$ is a network, or
\item\label{cone_gens_of_degree_2} $\omega$ has degree $2$, 
and satisfies the following three conditions determining $w_0$
  \begin{itemize}
    \item[]$e^*(\omega)=1$, for all cycle edges $e\in\ccE\setminus\ccL$
    \item[]$e^*(\omega)=2$, for an odd number of cycle legs,
    \item[]$e^*(\omega)=0$, for the remaining cycle legs.
  \end{itemize}
Each of the remaining components $\omega_i \in \tau(\ccG_i)$
is an element of degree at most two.
\end{enumerate}
\end{thm}

We postpone the proof until we prepare for it with some lemmas.
The idea of the proof is to use Lemma~\ref{cut_edge}
in order to remove the pendant trees $\ccG_1,\ldots,\ccG_k$
and work only with the polygon graph $\ccG_0$.
Lemma \ref{deg2points} describes all the degree~$2$ points
of the cone of a polygon graph and distinguishes the generators among them.

\begin{ex}\label{LMgens}
In Figure~\ref{pic_gens_of_LM} we illustrate the generators of the cone associated to the graph LittleMan
(one of the two graphs with one cycle and two leaves).
The first four are of degree $1$, the remaining three are of degree $2$.
\myVSPACEfigure\begin{center}
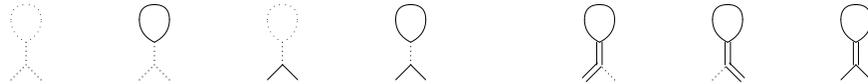
  \begin{minipage}{\textwidth}
\[
\begin{array}{cccc cccc cccc c}
  \LM....[][]
&\phantom{xx}&
  \LM-...[][]
&\phantom{xx}&
  \LM..--[][]
&\phantom{xx}&
  \LM-.--[][]
&\phantom{xxxxxx}&
  \LM-==.[][]
&\phantom{xx}&
  \LM-=.=[][]
&\phantom{xx}&
  \LM-=--[][]
\end{array}
\]  
\captionof{figure}{Generators of the cone of the LittleMan graph
\label{pic_gens_of_LM}
}
 \end{minipage} \end{center}\myVSPACEfigure
\end{ex}

\begin{lem}[Decomposition propagates to pendant trees]
\label{decomposition_propagates}
 Let $\ccG$ be any trivalent graph and $\omega \in \tau$ any cone element.
Let us also fix  a non-cycle edge $e \in \ccE \setminus \ccP$ 
such that $\ccG^e$ is a tree, where $\ccG^e=\ccG_1 \sqcup \ccG_2$.
Then any decomposition of $\rho^e_1(\omega)$ lifts to a decomposition of $\omega$.
\end{lem}
\begin{prf}
First note that both $\rho^e_1$ and $\rho^e_2$ preserve the degree, 
so  an element $\omega$ of degree $d$ in $\tau(\ccG)$
yields $\omega_1 \in \tau(\ccG_1)$ and $\omega_2 \in \tau(\ccG_2)$ both of degree $d$.
The semigroup of a tree is generated by networks, which are degree 1 elements, 
see Theorem~\ref{tree_gens}. 
This means that $\rho^e_2(\omega)$ is a sum of degree 1 elements.
Thus  if $\rho^e_1(\omega)$ can be decomposed,
then the same decomposition works for $\omega$ by choosing appropriate grouping 
of the summands of $\rho^e_2(\omega)$, because the degrees are preserved.
\end{prf}
\begin{cor}\label{polygon}
In the proof of Theorem~\ref{cone_gens} we can assume that 
the graph is a polygon graph.
\end{cor}
\begin{prf}
A graph with one cycle is a polygon graph with a tree  attached 
to each cycle leg $l_1,\ldots,l_k$.
We cut all the cycle legs  to obtain $k+1$ pieces of the graph $\ccG$:
a polygon graph $\ccG_0$ and $k$ trees: $\ccG_1,\ldots,\ccG_k$.
We denote by $\rho_0$ the composition of projections for each leg that we cut
$\rho_0 = \rho_0^{l_1} \circ \dotsb \circ \rho_0^{l_k}$.
Iteratively using Lemma~\ref{decomposition_propagates}
to decompose $\rho_0(\omega)$ we decompose~$\omega$.
\end{prf}
\begin{lem}[Degree $2$ elements of the cone]\label{deg2points}
Let $\ccG$ be a graph with exactly one cycle. 
Any degree $2$ element of $\tau(\ccG)$ 
except those in Theorem~\ref{cone_gens}~\ref{cone_gens_of_degree_2} 
is a sum of two networks.
\end{lem}

\begin{prf}
By Corollary~\ref{polygon} we may assume $\ccG$ is a polygon graph.
Let $\omega \in \tau(\ccG)$ a degree $2$ element.
The coefficient $e^*(\omega)$ of $\omega$ on any edge $e \in \ccE$ is $0,1$ or $2$.
We denote by
\[
\omega_{\ccE \setminus \ccP} := \{ e^*(\omega) \mid e \in \ccE \setminus \ccP \}
\]
the set of coefficients of $\omega$ on the cycle edges.
We distinguish between four types of $\omega$ based on $\omega_{\ccE \setminus \ccP}$.
For all but one we decompose $\omega$ as a sum of two networks.
  
If $0 \in \omega_{\ccE \setminus \ccP}$,  there exists a cycle edge 
$e \in \ccE(\ccG)$ with $e^*(\omega)=0$. We can cut it with no harm to $\omega$,
since $i^e(\omega) \in \tau(\ccG^e)$ is a degree $2$ element
in a cone of the trivalent tree $\ccG^e$, 
so it can be decomposed into a sum of degree $1$ elements.
This decomposition can be lifted to $\tau(\ccG)$, 
as we assumed $e^*(\omega)=0$. 
On Figure~\ref{pic_decomp_cycle_edge_zero} we show an example of 
this situation.

\myVSPACEfigure\begin{center}
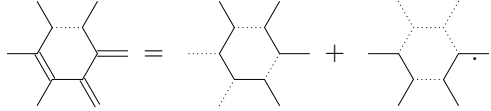
  \begin{minipage}{\textwidth}
\[
\xymatrix{
\begin{xy}
\HEXAGONnoxy{-.-=--}{=----=} 
\end{xy}
\quad =   
\begin{xy}
\HEXAGONnoxy{-.---.}{---..-} 
\end{xy}
\quad +
\begin{xy}
\HEXAGONnoxy{...-.-}{-..---} 
  \end{xy}.
}
\]
\captionof{figure}
{Decomposition of a degree two element with a weight zero cycle edge}
\label{pic_decomp_cycle_edge_zero}
\end{minipage} \end{center}\myVSPACEfigure
The second case is when $\omega_{\ccE\setminus \ccP} = \{2\}$, 
that is~$e^*(\omega)$ is~$2$ on all cycle edges. 
As $\omega$ has degree~$2$, we know that $\deg_v(\omega) \leq 2$ and as a consequence:
\[ e^*(\omega)=
  \begin{cases}
     2 &\text{if $e$ is a cycle edge}\\
     0 &\text{otherwise, i.e.~$e$ is a cycle leg.}
  \end{cases}
\]
Thus $\omega$ is twice the network consisting of all the cycle edges, 
as on the example on Figure~\ref{pic_decomp_cycle_edges_all_two}.

\myVSPACEfigure\begin{center}
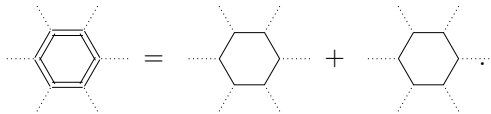
  \begin{minipage}{\textwidth}
\[
\xymatrix{
\HEXAGON{======}{......} \quad =
\HEXAGON{------}{......} \quad +
\HEXAGON{------}{......} 
}.
\]
\captionof{figure}
{Decomposition of a degree two element with all cycle edges of weight two}
\label{pic_decomp_cycle_edges_all_two}
\end{minipage} \end{center}\myVSPACEfigure

For the last two cases we let $l_1,\ldots,l_p$ be the set of all cycle legs 
with $l^*_i(\omega) \neq 0$ ordered anticlockwise and indexed 
by elements of $\ZZ_p$.  
Also we denote by $i \path j$ the path that starts at $l_i$ 
and goes anticlockwise along the intermediate cycle edges to $l_j$ 
and ends there.

In the third case  $\omega_{\ccE \setminus \ccP} = \{1\}$ and 
we will show that $\omega$ can be decomposed 
into a sum of two networks if and only if the number of cycle legs
for which $\omega$ has a coefficient $2$ is even. 
Suppose we have decomposed $\omega=\omega_1 + \omega_2$ into a sum of networks 
and neither $\omega_1$ nor $\omega_2$ contains the path consisting of all cycle edges 
(in which case the other one would be an empty network).
Then both $\omega_1$ and $\omega_2$ contain a positive even number of legs.
Any path in $\omega_1$ (respectively $\omega_2$) is of the type $i \path i+1$, 
from the cycle leg $l_i$ to the next one.
If the end was not the next one, there would be a path in $\omega_2$ 
(respectively $\omega_1$) containing an intermediate leg and as a result
there would be a cycle edge $e$ with $\omega(e)=2$, 
but $\omega_{\ccE \setminus \ccP}=\{1\}$.
Moreover, both $\omega_1$ and $\omega_2$ contain 
all non-zero cycle legs, each with value $1$, because 
they are both of degree $1$ and in the cone we have 
$\deg(\omega) \ge \deg_{\min}(\omega)$.
This is only possible when the number $p$ of non-zero cycle legs is even, 
and in that case we have the obvious decomposition with
$\omega_1=
\scriptsize{
\begin{array}{lclcl}
i_1 \path i_2  &+& \ldots &+& i_{p-1} \path i_{p} \\
\end{array}
}$
and 
$\omega_2=
\scriptsize{
\begin{array}{lclcl}
i_2 \path i_3  &+& \ldots &+& i_{p} \path i_{1} \\
\end{array}
}$.
Otherwise $\omega$ is a generator.
Examples of both these situations are drawn on 
Figure~\ref{pic_generator_and_cycle_edges_all_one}.

\myVSPACEfigure\begin{center}
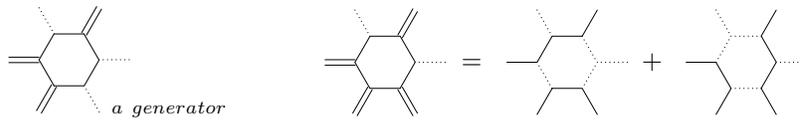
  \begin{minipage}{\textwidth}
\[
\xymatrix{
\HEXAGON{------}{.=.==.}_{a\ generator} &
\HEXAGON{------}{.=.===} \quad =
\HEXAGON{.--.-.}{.-.---} \quad +
\HEXAGON{-..-.-}{.-.---} 
}
\]
\captionof{figure}{Typical generator of degree two and 
a decomposition of a degree two element with all cycle edges of weight one}
\label{pic_generator_and_cycle_edges_all_one}
\end{minipage} \end{center}\myVSPACEfigure
In the last case $\omega_{\ccE \setminus \ccP} = \{1,2\}$.
When $l^*(\omega)=2$ we call $l$ a two-leg.
Denote by $l_{i_1},\ldots,l_{i_q}$ the subsequence of two-legs, 
numbered in such way that traveling along the cycle anticlockwise
from $l_{i_q}$ to $l_{i_1}$ there is a one-leg, provided that $\omega$ has some one-legs.
We observe that the number of one-legs between two consecutive 
two-legs is always even. 
This is best explained by drawing the picture from 
Figure~\ref{pic_deg_2_with_even_no_one_legs}.

\myVSPACEfigure\begin{center}
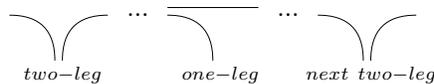
  \begin{minipage}{\textwidth}
\[
\begin{xy}
=(0, 0)"0";
=(-6, 0)"A";
=( 6, 0)"B";
=( 0,-6)"C";
=( 1,-8)"TAG";
(0,0);(1,0):; "TAG"*{\scriptstyle{two-leg}}; "A";"C" **\crv{"0"};  
(1,0);(2,0):; "C";"B" **\crv{"0"};
(10,0);(11,0):; "0"*{...};
(10,0);(11,0):; "TAG"*{\scriptstyle{one-leg}};
"A";"C" **\crv{"0"}; (0,1);(1,1):; "A";"B" **\dir{-};
(10,-1);(11,-1):; "0"*{...};
(10,0);(11,0):; "TAG"*{\scriptstyle{next\ two-leg}};
"A";"C" **\crv{"0"};  (1,0);(2,0):; "C";"B" **\crv{"0"};
\end{xy}
\]
\captionof{figure}{Element of degree two has even number of one-legs
between consecutive two-legs}
\label{pic_deg_2_with_even_no_one_legs}
\end{minipage} \end{center}\myVSPACEfigure
where the arcs are our $x_v,y_v,z_v$'s introduced in~\ref{local_paths}.
To produce an element of the lattice, 
the two on the same edge need to share the same value.
If there would be only zero-legs where the first dots are,
the local paths would not agree on some cycle edge.

The decomposition $\omega=(1,\omega_1)+(1,\omega_2)$ depends on 
the parity of $q$, which is the number of two-legs. 
We first work in the case $q=2r$ is even, the odd case 
uses the same idea with small modifications.
The Figure~\ref{pic_decomp_deg_2_even_two_legs}  
visualizes how the decomposition is constructed
in the case $q=4$.
\myVSPACEfigure\begin{center}
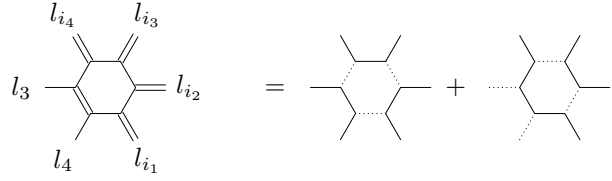
  \begin{minipage}{\textwidth}
\[
\xymatrix{
  \begin{xy}
\HEXAGONnoxy{---=--}{===--=}
\POS (0,0);(2.8,0):;
"A"*{l_{i_2}};
"B"*{l_{i_3}};
"C"*{l_{i_4}};
"D"*{l_3};
"E"*{l_4};
"F"*{l_{i_1}};
    \end{xy} &=
\HEXAGON{.-.-.-}{------} \quad +
\HEXAGON{-.---.}{---..-} 
}
\]
\captionof{figure}{Decomposition of a degree two element 
with even number of two-legs}
\label{pic_decomp_deg_2_even_two_legs}
\end{minipage} \end{center}\myVSPACEfigure
First we place all paths between two consecutive two-legs 
starting at an leg with an even (respectively odd) index into $\omega_1$ 
(respectively $\omega_2$).
Then, to take care of the one-legs, we add paths between consecutive one-legs 
lying between
$l_{i_{2j}}$ and $l_{i_{2j+1}}$ for some $j\in \{1,\ldots,p\}$ to $\omega_1$. 
Thus we get
\[
\omega_1=
\scriptsize{
\begin{array}{lclcl}
i_1 \path i_2 &+& (i_2+1) \path (i_2+2) &+ \ldots+& (i_3-2) \path (i_3-1)+ \\
i_3 \path i_4 &+& (i_4+1) \path (i_4+2) &+ \ldots+& (i_5-2) \path (i_5-1)+ \\
&&\vdots&&\\ 
i_{2r-1} \path i_{2r} &+& (i_{2r}+1) \path (i_{2r}+2) &+\ldots+& (i_1-2) \path (i_1-1)\\
\end{array}
}.\]
Similarly  we add paths between consecutive one-legs lying between
$l_{i_{2j-1}}$ and $l_{i_{2j}}$ for some $j\in \{1,\ldots,p\}$ to $\omega_2$. 
So we can write
\[
\omega_2=
\scriptsize{
\begin{array}{lclcl}
i_2 \path i_3 &+&  (i_3+1) \path (i_3+2) &+\ldots+& (i_4-2) \path (i_4-1)+\\
i_4 \path i_5 &+&  (i_5+1) \path (i_5+2) &+\ldots+& (i_6-2) \path (i_6-1)+\\
&&\vdots&&\\ 
i_{2r} \path i_1 &+& (i_1+1) \path (i_1+2) &+\ldots+& (i_2-2) \path (i_2-1)\\
\end{array}
}.\]
Clearly those paths in $\omega_1$ (resp.~$\omega_2$) are disjoint, so both are networks 
and by construction they yield a decomposition of our $\omega$.

When the number $q$ of two-legs is odd we need to 
adjust the above decomposition.
Again, we draw an example on 
Figure~\ref{pic_decomp_deg_2_odd_two_legs} for $q=3$
\myVSPACEfigure\begin{center}
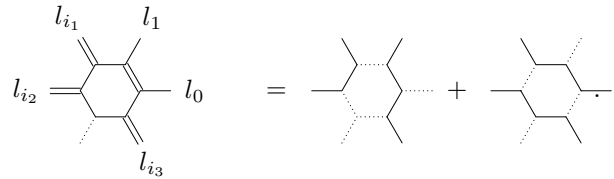
  \begin{minipage}{\textwidth}
\[
\xymatrix{
\begin{xy}
\HEXAGONnoxy{=-----}{--==.=}
\POS (0,0);(2.8,0):;
"A"*{l_0};
"B"*{l_1};
"C"*{l_{i_1}};
"D"*{l_{i_2}};
"E"*{};
"F"*{l_{i_3}};
\end{xy} &=
\HEXAGON{-.-..-}{.---.-} \quad +
\HEXAGON{--.--.}{-.--.-}.
}\]
\captionof{figure}{Decomposition of a degree two element with 
odd number of two-legs}
\label{pic_decomp_deg_2_odd_two_legs}
\end{minipage} \end{center}\myVSPACEfigure
As there is an odd number of two-legs, we will use two consecutive one-legs 
lying between $l_{i_q}$ and $l_{i_1}$ to make up for the missing two-leg, 
and we proceed as before to get
\[
\omega_1=
\scriptsize{
\begin{array}{lclcl}
i_1 \path i_2  &+& ( i_2+1) \path (i_2+2 ) +\ldots+ (i_3-2) \path (i_3-1)+\\
i_3 \path i_4  &+& (i_4+1) \path (i_4+2  ) +\ldots+ (i_5-2) \path (i_5-1)+\\
&& \vdots &&\\
i_{2r-1} \path i_{2r} &+& (i_{2r}+1) \path (i_{2r}+2) +\ldots+ (i_{2r+1}-2) \path (i_{2r+1}-1)\\
i_{2r+1} \path i_1-1&&&&\\
\end{array}
}\]
and
\[
\omega_2=
\scriptsize{
\begin{array}{lclcl}
i_2 \path i_3 &+& (i_3+1) \path (i_3+2) +\ldots+ (i_4-2) \path (i_4-1)+\\
i_4 \path i_5 &+& (i_5+1) \path (i_5+2) +\ldots+ (i_6-2) \path (i_6-1)+\\
&& \vdots &&\\
i_{2r} \path i_{2r+1} &+& (i_{2r+1}+1) \path (i_{2r+1}+2) +\ldots+ (i_1-4) \path (i_1-3)\\ 
i_1-2 \path i_1&&&&\\
\end{array}
}\]
This ends the proof of the lemma about the decomposable degree two 
elements of the cone.
\end{prf}

\begin{prf}[Proof of Theorem~\ref{cone_gens}]
The proof yields an algorithm for decomposing an arbitrary element $\omega$ 
of the cone $\tau(\ccG)$ into a sum of degree $1$ and $2$ generators.
By Corollary~\ref{polygon} we only need to prove the theorem 
when $\ccG$ is a polygon graph.
First we fix an orientation of the cycle of $\ccG$ and we call it 
anticlockwise in order to think of a planar embedding of the graph.
Let $\omega \in \tau(\ccG)$. We will find an element  $\mu$ 
of degree at most~$2$, such that $\omega-\mu \in \tau(\ccG)$.
Let  $v$ be a vertex and $l_v$ the cycle leg attached it.
We choose an map $i_v :\TRIPOD \ra \ccG$ so that 
the edge $e_3$ is mapped to $l_v$,
the edge $e_2$ is mapped to the edge  which points anticlockwise from $v$, and 
the edge $e_1$ is mapped to the edge  which points clockwise from $v$.
In this notation the coefficient of $\omega$ at $l_v$ is measured by $c_v$
and coefficient of the edge anticlockwise (respectively clockwise) from $v$ 
is measured by $a_v$ (respectively $b_v$). We will also use the local paths
$x_v, y_v, z_v$ defined in~\eqref{local_paths}.

Now we are ready to describe the algorithm to find $\mu$.
\begin{enumerate}

\item[Step 1.] \label{0edge} If there is a cycle edge $e$ with  $e^*(\omega)=0$,
we cut $e$ and obtain the graph $\ccG^e$ which is a trivalent tree.
Thus $\rho^e(\omega)$ is a sum of networks of path and this decomposition
lifts to graph $\ccG$.
\smallskip
\item[Step 2.] \label{1cycle} Otherwise  $e^*(\omega) \geq 1$ on all cycle edges $e$.
We set $\mu$ to have value $1$ on every cycle edge. 
Equivalently $\mu$ is defined by setting at every cycle vertex
$z_v(\mu)=1$, $x_v(\mu)=0$, $y_v(\mu)=0$.
Now if $\omega-\mu \in \tau(\ccG)$ we are done.
Otherwise $\omega-\mu$ fails one of the inequalities defining $\tau(\ccG)$.
It is not the one with degrees, since for each vertex $v\in\ccG$
we have $\deg_v(\mu) = 1$, which implies
$\deg_{\min}(\omega-\mu)  \leq \deg(\omega)-\deg(\mu) = \deg(\omega-\mu)$.
Clearly all coefficients of $\omega-\mu$ are positive.
Thus at some vertex $v \in \ccV$ our $\omega-\mu$
fails one of the triangle inequalities.
\smallskip
\item[Step 3.] \label{2legs} 
We will adjust $\mu$ to fix the triangle inequalities for $\omega-\mu$.
If a triangle inequality for $\omega-\mu$ at $v$ fails, 
then this is because $z_{v}(\omega)=0$.
In such a case we set $\mu(l_v)=2$, which will not make 
any coefficient of $\omega-\mu$ negative provided $c_v(\omega)\ge 2$.
But since $a_{v}(\omega), b_{v}(\omega) \ge 1$ and $z_{v}(\omega)=0$ 
we must have  $x_{v}(\omega), y_{v}(\omega) \ge 1$.
This implies $c_v(\omega)\ge 2$ as required.
In terms of $x_{v},y_{v},z_{v}$ we have
decreased $z_{v}(\mu)$ by one and 
increased both  $x_{v}(\mu)$ and $y_{v}(\mu)$ by one.
\smallskip
\item[Step 4.] \label{deg} 
We need a last adjustment on  $\mu$  
to assure the additivity of degree where it is attained, i.e.
for any $v$ such that $\deg_v(\omega)=\deg(\omega)$ 
we need $\deg_v(\mu)=2$ since $\deg(\mu)=2$.
This is to  ensure $\omega-\mu \in \tau$.
We call $v$ \emph{degree deficient vertex} if $\deg_v(\omega)=\deg(\omega)$ 
and  $\deg_v(\mu)=1$.

If $v$ is degree deficient and in addition $x_v(\omega)>0$ and $y_v(\omega)>0$ we set 
$x_v(\mu)=y_v(\mu)=1$ and $z_v(\mu)=0$.

If $v$ is degree deficient and both $x_v(\omega)=y_v(\omega)=0$ are zero, 
then $z_v(\omega) = \deg_v(\omega)$, so at both next and previous 
cycle vertex the degree is attained
\begin{multline}\label{deg_deficient}
\deg_{v_{next}}(\omega)=\deg_{v_{prev}}(\omega)=\deg(\omega)\\
\text{since }\qquad z_{v_{next}}+y_{v_{next}}=z_v+x_v=z_{v_{prev}}+x_{v_{prev}}
\end{multline}
If all degree deficient  vertices were of this type, 
then $\omega$ was a multiple of the path consisting of all cycle edges.

Now we divide the set of all deficient vertices (which all have
at least one of $x_v(\omega)$ or $y_v(\omega)$ equal to zero)
into  sequences of adjacent ones.
Let us fix our attention to such a sequence 
(we already excluded the case when it has the same end and beginning).
Call it  $v_1,\ldots,v_r$. Then because all $v_i$'s are degree deficient:
$\deg(\omega)=\deg_{v_i}(\omega)$. The last one has $y_{v_r}(\omega)>0$ 
by \eqref{deg_deficient}
\[
\deg(\omega)=z_{v_{r-1}}(\omega)=y_{v_r}(\omega)+z_{v_r}(\omega) \leq z_{v_r}(\omega) < \deg(\omega)
\]
In the same way $x_{v_1}(\omega)>0$. Finally all the middle ones have $z_{v_i}(\omega)>1$.
This implies that we can redefine $\mu$ on our fixed sequence $v_1,\ldots,v_r$ 
preserving all other properties and fixing the degree deficiency:
\begin{align*}
x_{v_1}(\mu)&=1&
x_{v_i}(\mu)&=0&
x_{v_r}(\mu)&=0\\
y_{v_1}(\mu)&=0&
y_{v_i}(\mu)&=0&
y_{v_r}(\mu)&=1\\
z_{v_1}(\mu)&=1&
z_{v_i}(\mu)&=2& 
z_{v_r}(\mu)&=1\\
\end{align*}
where $i\in \{2,\ldots,r-1\}$.
We do this for all such sequences and we have the required~$\mu$.
\end{enumerate}
Now $\mu$ is an element of degree $2$ and is either a generator 
or can be decomposed into a sum of two generators in degree 1, 
as described in Lemma~\ref{deg2points}.
\end{prf}

\begin{rmk}
If we allow more cycles, we can have generators of higher degree.
As we can see on Figure~\ref{pic_ex_deg_3_gen} 
the graph with two loops and one leaf has a degree three generator:
one on the two loops, two on the three other edges.

\myVSPACEfigure\begin{center}
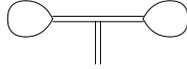
  \begin{minipage}{\textwidth}
\[
\begin{xy}
(0,0);(0,0.6):;0*{\CHERRY};
(0,0);(0,-20) **\dir{=};
(-0.5,-10);(-10,-10) **\dir{=};
(0,-20);(-1,-20):;0*{\CHERRY};
\end{xy}
\]
\captionof{figure}{Example of a degree three generator}
\label{pic_ex_deg_3_gen}
\end{minipage} \end{center}\myVSPACEfigure

\end{rmk}

\subsection{Embedding.}\label{section_embedding}

The aim of this section is to find a common ambient space
for phylogenetic models of all graphs with the same topological invariants.
The way to construct this ambient space
follows easily form the definition of the model.

\begin{thm}\label{embedding}
The phylogenetic model of a trivalent graph $\ccG$ embeds in 
a projective toric variety $\PP_{g,n}$, 
which is a good quotient of projective space by 
an action of a $g$-dimensional torus. 
This action as well as the variety  $\PP_{g,n}$
depends only on the first Betti number 
and number of leaves of  $\ccG$, up to reordering of coordinates.
\end{thm}
\begin{prf}
Models of trivalent trees with $N$ leaves embed naturally in $\PP^{2^{N-1}-1}$, 
with coordinates $x_\kappa$ where $\kappa \in \{0,1\}^N$ is a sequence of length $N$ 
with values in $\{0,1\}$ and even number of $1$-entries, 
see Proposition~\ref{tree_gens}.
We cut $g$ cycle edges of the graph $\ccG$ so that we obtain
a trivalent tree $\ccT$ with set $S$ of $M=n+2g$ leaves.
Dividing $\PP^{2^{M-1}-1}$ by the action of the $g$-dimensional torus $\TT(g,n)$
that corresponds to gluing these leaves back together,
yields the required embeddings.

We will now explain that this action does not depend 
on the graph~$\ccG$, up to choice of coordinates on 
Let us label the set of leaves  by $1\ldots M$ 
and divide it into three disjoint sets
$S=S_0 \sqcup S_1 \sqcup S_{-1}$ as follows.
$S_0$ is the set of leaves of the original graph $\ccG$.
The $2g$ new leaves of $\ccT$ come in pairs $(l, l^-)$, 
where both $l$ and $l^-$  used to be the same edge in $\ccG$.
We put $l$ in $S_1$ and $l^-$ in  $S_{-1}$. 

The action of the $g$-dimensional torus $\TT(g,n)$ on  $\PP^{2^{M-1}-1}$ 
is given by a matrix 
\[
\{ \zeta^l_{\kappa}\}^{l=1 \ldots g}_{\kappa \in \{0,1\}^M} \text{ , where }\quad \zeta^l_\kappa = \kappa_l-\kappa_{l^-} 
\]
 Thus this action is independent of the shape of $\ccG$ and
depends only on $(g,n)$ up to choice of order of coordinates.
\end{prf}
We will now illustrate Theorem~\ref{embedding}.
\begin{ex}\label{LMembedding}
We describe $X \left ( \LMfont \right )$, the model of LittleMan, 
together with its embedding into $\PP_{1,2}$.
We know that   the model of a trivalent tree $\TREEfont$ 
is a complete intersection of two quadrics, 
see \cite[Ex.~2.6]{buczynska_wisniewski},
\begin{multline}\label{tree_equations}
(x_{1100}x_{0011} - x_{0000}x_{1111}, \quad
x_{1001}x_{0110} - x_{1010}x_{0101} )  \\
\subset
\Proj \CC 
[x_{0000},x_{1100},x_{0011},x_{1111},
 x_{1010},x_{1001},x_{0101},x_{0110}]
\end{multline}
where each variable $x_\kappa$ corresponds to a degree $1$ generator $\tau(\TREEfont)$.
By Proposition~\ref{tree_gens} a generator is given by its values on the leaves, 
which form the index $\kappa \in \{0,1\}^4$.

To obtain $X \left ( \LMfont \right )$, we glue two leaves, 
say leaf $e_1$ is glued with the leaf $e_2$ to get the loop in LittleMan
as shown on Figure~\ref{pic_gluing_LM}.

\myVSPACEfigure\begin{center}  \begin{minipage}{\textwidth}
\[\xymatrix{
\TREE-----[1][2][3][4]
\ar@{~>}[r]
&
\LM----[3][4]
}\]
\captionof{figure}{Obtaining LittleMan graph from a four-leaf tree}
\label{pic_gluing_LM}
\end{minipage} \end{center}\myVSPACEfigure

Now $X \left ( \LMfont \right )$ is a GIT-quotient of $X \left ( \TREEfont \right )$
by the action $\lambda^{e_1-e_2}$. 

In the same way it embeds in a GIT-quotient of  $\PP^7$. One easily sees that
\[
\PP^7  \git  \CC^* = 
V(y_1 \cdot y_2 - z_1 \cdot z_2) \subset
\PP(1^4, 2^4)
\]
where 
\[
\PP(1^4, 2^4)=
\Proj \CC 
[x_{0000},x_{1100},x_{0011},x_{1111},
y_1, y_2, z_1, z_2] 
\]
is a weighted projective space and
\[
y_1=x_{1001}\cdot x_{0110}, \quad y_2=x_{1010}\cdot x_{0101}, \quad
z_1=x_{1010}\cdot x_{0110}, \quad z_2=x_{0101}\cdot x_{1001}
\]
are the $\CC^*$--invariant variables of degree $2$. From Theorem~\ref{GIT} it follows that $X \left ( \LMfont \right )$ is given by Equations~\eqref{tree_equations} 
of $X \left ( \TREEfont \right )$ in $\PP^7$, 
rewritten in the coordinates of $\PP^7  \git  \CC^*$.
The second one becomes $y_1-y_2=0$, so
\[\begin{array}{rrcrcrc}
X \left ( \LMfont \right ) &=& 
\Proj \CC[\tau \left ( \LMfont \right )] &=& 
V(x_{1100}x_{0011} - x_{0000}x_{1111}, y_1^2-z_1z_2) 
&\subset& \PP(1^4,2^3) \\
&&&&&&\parallel\\
&&\bigcap&&&&V(y_1-y_2)\\
&&&&&&\bigcap\\
\PP_{1,2}&=&\PP^7  \git  \CC^* &=&
V (y_1y_2-z_1z_2) 
&\subset&\PP(1^4, 2^4) \\
\end{array}\]

If we replace each variable in the equations of $X(\LMfont)$
by its representation on the graph, we get the picture
from Figure~\ref{pic_ideal_LM}.
\myVSPACEfigure\begin{center}  \begin{minipage}{\textwidth}
\[
\LM-...[][] + \LM..--[][] = \LM....[][] + \LM-.--[][] 
\qquad \qquad
\LM-=--[][] + \LM-=--[][] = \LM-==.[][] + \LM-=.=[][] 
\]
\captionof{figure}{Ideal of the model LittleMan graph}
\label{pic_ideal_LM}
\end{minipage} \end{center}\myVSPACEfigure
\end{ex}

\begin{ex}\label{HAMMOCKembedding}
As in the previous example, we work out properties of 
$X \left ( \HAMMOCKfont \right )$ --- the model of Hammock.
We use the same action $\lambda^{l_1-l_2}$, but we change the embedding
$X \left ( \TREEfont \right ) \mono \PP^7$ by relabeling the leaves 
as on Figure~\ref{pic_gluing_Hammock}.

\myVSPACEfigure\begin{center}  \begin{minipage}{\textwidth}
\[\xymatrix{
\TREE-----[1][3][2][4]
\ar@{~>}[r]
&\quad
\HAMMOCK----[4][3]
}\]
\captionof{figure}{Obtaining Hammock graph from a four-leaf tree}
\label{pic_gluing_Hammock}
\end{minipage} \end{center}\myVSPACEfigure

As the labeling of $\TREEfont$ was modified,
Equations~\eqref{tree_equations} become
\[
V(x_{1010}x_{0101} - x_{0000}x_{1111}, \quad
 x_{1001}x_{0110} - x_{1100}x_{0011}  )  
\subset \PP^7.
\]
We again rewrite them in the invariant coordinates of $\PP(1^4,2^4)$ to get:

\[\begin{array}{rrcrcrc}
X \left ( \HAMMOCKfont \right ) &=& 
\Proj \CC[\tau \left ( \HAMMOCKfont \right)] &=&
V( x_{0000}x_{1111}x_{1100}x_{0011}-z_1z_2)  
&\subset&  \PP(1^4,2^2) \\
&&&&&&\parallel\\
&&\bigcap&&
\multicolumn{3}{r}{V(y_1 - x_{0000}x_{1111},  x_{1100}x_{0011} - y_2 )}\\
&&&&&&\bigcap\\
\PP_{1,2}&=&\PP^7  \git  \CC^* &=&
V (y_1y_2-z_1z_2) 
&\subset&\PP(1^4, 2^4) \\
\end{array}\]

If we replace each variable in the degree $4$ equation of $X(\HAMMOCKfont)$
by its representation on the graph, we get the picture
shown on Figure~\ref{pic_ideal_Hammock}.

\myVSPACEfigure\begin{center}  \begin{minipage}{\textwidth}
\[\xymatrix{
\HAMMOCK....[][]
\quad+\quad
\HAMMOCK.--.[][]
\quad+\quad
\HAMMOCK--.-[][]
\quad+\quad
\HAMMOCK-.--[][]
\quad=\quad
\HAMMOCK=--.[][]
\quad+\quad
\HAMMOCK.--=[][]
}\]
\captionof{figure}{Ideal of the model of the Hammock graph}
\label{pic_ideal_Hammock}
\end{minipage} \end{center}\myVSPACEfigure
\end{ex}

\section{Flat families}\label{section_defomation}

Models of trivalent trees that differ by one mutation 
live in a flat family in a projective space \cite{buczynska_wisniewski}. This statement
almost remains true for trivalent graphs, by the same argument. The only
difference is that we get a family in the projective toric variety $\PP_{g,n}$ 
instead of a usual projective space.

Recall that in Remark~\ref{torus_of_leaves} we have associated 
to a subset of leaves $\ccS\subset \ccL$ with $k$ elements 
an action of the $k$-dimensional torus $\TT(\ccS)$ on the model $X(\ccG)$.
\subsection{Key examples.}\label{section_deformation_examples}
We construct two-dimensional flat families containing models of small graphs.
They become the building blocks for deformations of bigger graphs.

\begin{ex}[g=0, n=4]\label{tree_family} Let $ \TREEfont$ be a trivalent tree
with four leaves. 
In  \cite[Ex. 2.20]{buczynska_wisniewski} we constructed a flat family
\[
\ccX^0 \mono \ccB \times \PP^7,
\] 
where  
\begin{itemize}

\item $\ccB$ is an open subset of $\PP^2$ with coordinates 
$b_{(1,2)(3,4)}$, $b_{(1,3)(2,4)}$, $b_{(1,4)(2,3)}$

\item the torus $\TT(\ccL)$ acts on $\ccB \times \PP^7$ via the second coordinate,
that is, for a leaf $l$ of $\TREEfont$ and coordinate $x_{\kappa}$ we have
$\lambda_{v_i}(t)(b_{(.)(.)},x_\kappa)=(b_{(.)(.)},t^{\kappa(l)}x_\kappa)$,
\item the equidimentional projection \mbox{$\ccX^0\ra\ccB$} 
contains the three special fibers 
$\ccX^0_{[1,0,0]}$, $\ccX^0_{[0,1,0]}$ and $\ccX^0_{[0,0,1]}$
which  are models of aforementioned four-leaf trees,

\item $\ccX^0$ is a $\TT(\ccL)$-invariant complete intersection in  
$\ccB\times\PP^7$ of the two quadrics
\vskip -1cm
\begin{multline*}
b_{(12)(34)}\cdot x_{1100}x_{0011} \,+\,
b_{(13)(24)}\cdot x_{1010}x_{0101} \,+\,
b_{(14)(23)}\cdot x_{1001}x_{0110}\\
\shoveright{\,-\,\left(b_{(12)(34)}+b_{(13)(24)}+b_{(14)(23)}\right) \cdot x_{0000}x_{1111}}\\
\phantom{haha}\\
\shoveleft{\left(b_{(13)(24)}-b_{(14)(23)}\right)\cdot x_{1100}x_{0011}
+\left(b_{(14)(23)}-b_{(12)(34)}\right)\cdot x_{1010}x_{0101}}\\
+\left(b_{(12)(34)}-b_{(13)(24)}\right)\cdot x_{1001}x_{0110}.
\end{multline*}
\end{itemize}
\end{ex}
\begin{ex}[$g=1$, $n=2$]\label{LM_family}
We construct a family which contains models of graphs with one cycle and two leaves.
It arises as a $\CC^*$-quotient of the family  $\ccX^0 \mono \PP^7$  
from Example~\ref{tree_family}.
Let us fix a $\CC^*$-action  $\lambda^{l_1-l_2}$ on the ambient $\PP^7$, 
thus on $\ccX^0$ and tree models as well, by choosing leaves 
$l_1$ and $l_2$ labeled by $1$ and $2$ respectively.
Each of the three trees yields a graph, when two leaves are glued together.
Up to graph isomorphism, there are two graphs with one cycle and two leaves.
As we are mutating along fixed edge, the LittleMan appears once
and the Hammock twice.
The picture of the three possible trees becomes
the one on Figure~\ref{pic_mutations_g1_2l}.

\myVSPACEfigure\begin{center}
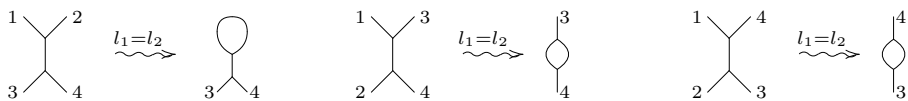
  \begin{minipage}{\textwidth}
\[ \xymatrix{
\TREE-----[1][2][3][4]
\ar@{~>}[r]^{l_1=l_2}
\quad & \quad
\LM----[3][4] 
&\quad
\TREE-----[1][3][2][4]
{\ar^{\ \ \ l_1=l_2}@{~>}[r]}
\quad & \quad
\HAMMOCK----[4][3]
&\qquad
\TREE-----[1][4][2][3]
{\ar@{~>}[r]^{\ \ \ \ \ l_1=l_2}}
\quad & \quad
\HAMMOCK----[3][4]
}\]
\captionof{figure}{Mutations of graphs with one cycle and two leaves}
\label{pic_mutations_g1_2l}
\end{minipage} \end{center}\myVSPACEfigure
The new ambient space 
$\PP^7  \git  \CC^* = (y_1\cdot y_2 - z_1\cdot z_2) \subset \PP(1^4, 2^4)$
was discussed in Example~\ref{LMembedding}. 

By Theorem~\ref{GIT}, the new family 
$\ccX^0  \git  \CC^*$ is given by the same, 
$\TT(S)$-invariant equations of $\ccX^0$. 
We rewrite them in the (invariant) coordinates of $\PP(1^4,2^4)$
\begin{multline*}
b_{(12)(34)}\cdot x_{1100}x_{0011} +b_{(13)(24)}\cdot y_2+b_{(14)(23)}\cdot y_1
-(b_{(12)(23)}+b_{(13)(24)}+b_{(14)(23)})\cdot x_{0000}x_{1111},\\ 
\shoveleft{(b_{(23)(14)}-b_{(14)(23)})\cdot x_{1100}x_{0011}
+(b_{(14)(23)}-b_{(12)(34)})\cdot y_2  
+(b_{(12)(34)}-b_{(23)(14)})\cdot y_1,} \\
\shoveleft{y_1\cdot y_2-z_1\cdot z_2.}\\
\end{multline*}
\vskip -0.5cm
To understand how this works a little better, 
let us look at particular coordinate of $\PP(1^4,2^4)$, 
say $y_2=x_{1001}x_{0110}$,
and draw on Figure~\ref{pic_coordinates_under_mutation}
its representation for each graph.
\myVSPACEfigure\begin{center}
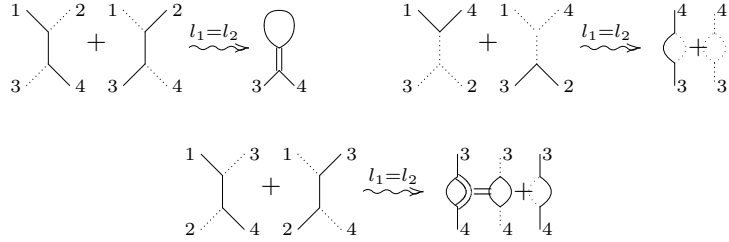
  \begin{minipage}{\textwidth}
\begin{gather*}
\xymatrix{
\TREE-.-.-[1][2][3][4] 
+
\TREE.---.[1][2][3][4] 
\ar@{~>}[r]^{\ \ \ \ \ \ l_1=l_2}
&
\LM-=--[3][4] 
&\quad
\TREE--...[1][4][3][2]
+
\TREE...--[1][4][3][2]
\ar@{~>}[r]^{\ \ \ \  \ \ l_1=l_2}
& \quad
\HAMMOCK--.-[3][4]
+
\HAMMOCK....[3][4]
}  \\ \phantom{empty line}\\ \xymatrix{ 
\TREE-.-.-[1][3][2][4] 
+
\TREE.---.[1][3][2][4] 
\ar@{~>}[r]^{\ \ \ l_1=l_2}
&\quad
\HAMMOCK--=-[4][3]
=
\HAMMOCK.--.[4][3]
+
\HAMMOCK-.--[4][3]
}\end{gather*}
\captionof{figure}{Behavior of coordinates of $\PP_{1,2}$ under mutations}
\label{pic_coordinates_under_mutation}
\end{minipage} \end{center}\myVSPACEfigure
\end{ex}

\subsection{General case.}\label{flat_ci_family}
We construct a flat family containing models of 
all mutations of a  given trivalent graph $\ccG$ 
along a fixed edge $e \in \ccE \setminus \ccP$. 
We follow \cite[Const.~3.5]{buczynska_wisniewski} replacing 
polytopes by cones.

Let $\ccG$ be a graph with an inner edge $e_0$ which 
contains two trivalent inner vertices.  
We can write $\ccG$ as a sum  of 
a not necessarily connected graph $\ccG_1$  with a set $S$ of $k$  
distinguished leaves $l_1,\ldots,l_k$ with $k\in\{0,\ldots,4\}$
and a graph $\ccG_0$ having the edge $e_0$ as its inner edge
and $\ccG_0$ is one of the following three small graphs

\begin{enumerate}

\item a tree $\TREEfont$ with four leaves $v_1,\ldots,v_4$, 
where the edge incident to the leaf $l_i$ is identified 
with the edge incident to the leaf $v_i$,

\item LittleMan $\LMfont$ or Hammock $\HAMMOCKfont$ --- a graph with four edges 
and two leaves $v_1$ and $v_2$,

\item Dumbell $\DUMBBELLfont$ a graph with three edges 
and two loops and no leaves.

\end{enumerate}
From Proposition~\ref{cone_of_graft_is_fiberproduct} 
we can compute the lattice $\Mg(\ccG)$ and the cone $\tau(\ccG)$ 
from those of the pieces $\ccG_0$ and $\ccG_1$
\begin{align*}
\Mg(\ccG)&=\Mg(\ccG_0)\times_g\Mg(\ccG_1)\cap
\bigcap_{i=1}^k\ker(\ell_i-v_i)\\
\tau(\ccG)&=\tau(\ccG_0)\times_g\tau(\ccG_1) \cap
\bigcap_{i=1}^k\ker(\ell_i-v_i). \\
\end{align*}

We consider the lattice $\Mg_{amb}$ and the cone $\tau_{amb}$, which are 
the combinatorial data of the ambient space $\PP_{g(\ccG_0), n(\ccG_0)}$ 
described by Theorem~\ref{embedding},
Example~\ref{LMembedding} and Example~\ref{HAMMOCKembedding}.
The inclusion $X(\ccG_0) \subset \PP_{g(\ccG_0), n(\ccG_0)}$ yields maps:
$\Mg_{amb} \ra \Mg(\ccG_0)$ and $\tau_{amb} \ra \tau(\ccG_0)$.
Forms $v_i$, $i=1,\dots k$ pull-back to $\Mg_{amb}$ and we
denote them by $\tilde{v}_i$, respectively. Now we define

\begin{align*}
\Mg_{\ccY} &=\Mg_{amb}(\ccG_0)\times_g \Mg(\ccG_1) \cap 
\bigcap_{i=1}^k\ker(\ell_i-\tilde{v}_i)  \\
\tau_{\ccY} &=\tau_{amb} (\ccG_0)\times_g\tau(\ccG_1) \cap 
\bigcap_{i=1}^k\ker(\ell_i-\tilde{v}_i). \\
\end{align*}
  
We define a toric variety $\ccY=\Proj \CC[\tau_{\ccY}]$. 
Since the good quotient is a categorical quotient,
by the construction we have the embedding
\[X(\ccG)\hookrightarrow \ccY.\]

\begin{lem}\label{quotient-of-product}
The inclusions 
\[
\begin{array}{ccc}
\Mg_{\ccY} \hookrightarrow 
\Mg_{amb} \times\ \Mg(\ccG_1)&
{\rm and}&
\tau_{\ccY}\hookrightarrow
\tau_{amb} \times \tau(\ccG_1) 
\end{array}
\]
induce a rational map 
\[
\PP_{g(\ccG_0),n(\ccG_0)}\times X(\ccG_1) \dashrightarrow \ccY
\]
which is a good quotient map (of the set over which it is defined) with respect
to the action of the k-dimensional torus $\TT_0$ generated by
one-parameter groups $\lambda_{v_i-\ell_i}$, where $i=1,\dots k$. The
subvariety
$$\widehat{\ccX}=\ccX^0\times X(\ccG_1)\hookrightarrow \ccB\times\PP^7\times X(\ccG_1)$$ 
is $\TT_0$-invariant and its quotient $\ccX$ is
locally complete intersection in $\ccB\times\ccY$.
\end{lem}
\begin{prf}
The map given by inclusions of cones and lattices is 
a good quotient map by Theorem~\ref{toricGIT}. 
Invariance of the resulting subvariety $\widehat{\ccX}$ 
follows by the invariance of $\ccX^0\hookrightarrow\ccB\times\PP^7$ 
discussed in Example~\ref{tree_family}. 
Finally, since $\widehat{\ccX}$ is a complete
intersection in $\ccB\times\PP^7\times X(\ccG_1)$
its image $\ccX$ is a locally complete intersection in the quotient $\ccB\times\ccY$.
This follows from the definition of good
quotient, which locally is an affine quotient, \cite[Ch. 5]{bialynicki_enc},
hence functions defining $\widehat{\ccX}$ locally descend to
functions defining $\ccX$.
\end{prf}

\begin{lem}\label{flatness-of-family}
Over an open set $\ccB'\subset\PP^2$ containing points $[1,0,0]$,
$[0,1,0]$, $[0,0,1]$ the projection morphism $\ccX\ra\ccB'$ is flat. The
fibers over points $[1,0,0]$, $[0,1,0]$, $[0,0,1]$ are reduced and
isomorphic to, respectively, the geometric model of $\ccG$ and of its
elementary mutations along the edge $e_0$.
\end{lem}
\begin{prf}
First we note that the fibers in question, $\ccX_{[*,*,*]}$, of
$\ccX\ra\ccB$ are geometric models as we claim. Indeed this follows from
the universal properties of good quotients, c.f.~\cite{bialynicki_enc}, as they
are quotients of the respective products
$\ccX^0_{[*,*,*]} \times  X(\ccG_1)$, which are located, as
three invariant subvarieties, in
$\widehat\ccX=\ccX^0\times  X(\ccG_1)$.  This, in
particular, implies that the respective fibers of $\ccX\ra\ccB$ are of
the expected dimension, hence they are contained in a set
$\ccB'\subset\PP^2$ over which the map in question is equidimentional.
Since $\ccY$ is toric it is Cohen-Macaulay and because $\ccX$ is a locally
complete intersection in $\ccY$, it is Cohen-Macaulay too
\cite[Prop.~18.13]{eisenbud}. Finally, the map $\ccX\ra\ccB'$ is
equidimentional hence it is flat, because $\ccB'$ is smooth, see
\cite[Thm.~18.16]{eisenbud}
\end{prf}

\begin{thm}\label{deformation_equivalent}
Geometric models of connected trivalent graphs with the same number 
of leaves $n$ and cycles $g$ are deformation equivalent in 
the projective toric variety $\PP_{g,n}$, which is a quotient of 
$\PP^{2^{n+2g-1}-1}$ by a $g$-dimensional torus.
\end{thm}

\begin{prf}
This is a combination of Proposition~\ref{graphs_mutation_equivalent} and of
Lemma~\ref{flatness-of-family}.
\end{prf}

\begin{rmk}
Since the phylogenetic model of disjoint union of graphs is a product
of the models of the pieces, see~\ref{cone_of_graft_is_fiberproduct},
we have proved that models of graphs with the same topological invariants
are deformation equivalent.
\end{rmk}
\section{Hilbert function of the model.}\label{section_Hilber_function}

As we mentioned in Section~\ref{section_defomation}, 
given a projective variety embedded in a projective space 
we have Hilbert function coming from the grading of the coordinates ring,
or equivalently from the action of the one-dimensional torus $\CC^*$.
If our projective variety is equipped with an action of a bigger 
torus it is natural 
(see \cite{haiman_sturmfels_multigraded_hilbert_schemes}) 
to consider a \emph{multigraded Hilbert function}, 
whose domain consists of the characters of the torus.

We study multigraded Hilbert function of an embedded projective toric variety
with the multigrading given by a subtorus of the big torus.
For the graph model $X(\ccG)$  the subtorus comes from a subset of leaves.
We prove in Theorem~\ref{same_hs} that the Hilbert function depends only 
on the topological invariants of the graph by showing 
that deformations constructed in 
the proof of Theorem~\ref{deformation_equivalent} 
preserve the whole Hilbert series. 
As tools we use Lemma~\ref{hs_of_ci} to compute the Hilbert series 
of a torus invariant complete intersection.
Lemma~\ref{hs_of_git} gives the formula for the Hilbert series 
of a quotient of a toric variety by a subtorus of the big torus.
We first state them in the algebraic setting.

\subsection{Ring with a torus action.}\label{section_rings_with_torus_action}
Let $R$ be a commutative $\CC$-algebra with an action of a torus $\TT$.
Let $M_{\TT}=\{\chi : \TT \ra \CC^*\}$
denote the group of characters of the torus $\TT$.
Then we can write 
\[
R = \bigoplus_{\chi \in M_{\TT} } R_{\chi}
\]
as a sum of isotypical pieces indexed by the characters of the torus. 
We assume that each $R_{\chi}$ has finite dimension over $\CC$.  
Then its \textbf{Hilbert function} $H_{R,\TT} : M_{\TT} \ra \NN$ is 
\[
H_{R,\TT}(\chi) :=  \dim R_{\chi}. 
\]
The \textbf{Hilbert series} is the generating series of $H_{R,\TT}$
\[
h_{R,\TT} (t):= \sum_{\chi \in M_{\TT}} \dim R_{\chi} \cdot t^{\chi}.
\]

\begin{lem}\label{hs_of_regular_sequence}
Let $R = \bigoplus_{\chi \in M_{\TT}} R_\chi$ be a ring with a torus action.
If $f_1,\ldots,f_q$ are homogeneous with $f_i \in R_{\chi_i}$
and form a regular sequence in $R$, then
\[
h_{R / \langle f_1,\ldots, f_q\rangle, \TT} (t)= h_R(t) \cdot (1-t^{\chi_1})\cdot \ldots \cdot(1-t^{\chi_q}).
\]
\end{lem}
\begin{prf}
The statement for the single graded Hilbert series is explicitly given 
in \cite[Cor.~3.2]{stanley}. Its multigraded, more general,
with minor additional assumption can be found
in \cite[Claim 13.38]{sturmfels_miller}.
The lemma follows by induction on the length of the regular sequence.
For any homogenous $f \in R_{\chi_f}$  which is not  
a zero divisor in $R$ and any $\chi \in M$ we have the exact sequence
of $\TT$-modules
\[ 
\xymatrix{
0 \ar[r] &  R_{\chi - \chi_f} \ar[r]^-{f \cdot} & 
R_\chi  \ar[r] &\left( R / (f) \right)_\chi \ar[r] & 0,
} 
\]
which implies that
\[ 
H_{R/(f)}(\chi) = 
\dim \left( R / (f) \right)_\chi = \dim R_\chi - \dim R_{\chi - \chi_f}
= H_R(\chi) - H_R(\chi - \chi_f)
\]
This is equivalent to the required equality for Hilbert series.
\end{prf}

Given a subtorus $\iota : \boldS \mono \TT$ we have the corresponding epimorphism 
of the character groups $\iota^* : M_{\TT} \epi M_{\boldS}$ 
and we can form the $\boldS$-invariant subring $R^{\boldS}$ of $R$
equipped with the residual action of the quotient torus $\TT/\boldS$
\[
R^{\boldS}=\bigoplus_{\chi \in M_{\TT/\boldS}} R_{\chi}.
\]
Then we have the following formula for the Hilbert series.
\begin{lem}\label{hs_of_invariants}
Let $R$ be a ring equipped with an action of a torus $\TT$ 
and let $\iota: \SS \mono \TT$ be a subtorus.
Then the Hilbert series of the invariant ring $R^{\SS}$ is
\[
h_{R^{\SS},\TT/\SS} = \sum_{\chi \in \ker \iota^* } t^{\chi}\dim R_{\chi},
\]
where $\iota^* : M_{\TT} \epi M_{\SS}$ is the dual map of the character groups.
\end{lem}

\subsection{Toric variety with a distinguished subtorus.}
\label{section_toric_with_subtorus_action}
Let $X = \Proj R$ be a projective toric variety of dimension~$d$ 
with an ample Weil divisor~$\ccL$ where 
$R = \bigoplus_{m \in \NN} R_m = \bigoplus_{m \in \NN} H^0(X, m\ccL)$
as in Section~\ref{section_GIT}.
Then $R$ has an action of a $d+1$ dimensional torus which is the product of
the $d$-dimentional torus $\TT$ of $X$ and the $\CC^*$ from the grading.
Any subtorus $\SS \mono \TT$ of dimension $r$ induces a $\ZZ^{r+1}$-sub-grading.
Then its \textbf{multigraded Hilbert function} 
$H_{X,\SS} : M_{\CC^* \times \SS} \ra \NN$ \textbf{with respect to $\SS$} is 
\[
H_{X, \SS}(\chi) := H_{R,\SS} (\chi) = \dim R_{\chi}.
\]
The generating series of $h_{R,\SS}$ is the 
\textbf{multigraded Hilbert series with respect to~{$\SS$}}
\[
h_{X,\SS} (t) := h_R (t) = \sum_{\chi \in M_{\SS}} \dim R_{\chi} \cdot t^{\chi}.
\]
We have the  following corollary of Lemma~\ref{hs_of_regular_sequence}.

\begin{cor}\label{hs_of_ci}
Let $Y=\Proj R$ be a projective toric variety with an action of 
an $r$-dimensional subtorus  $\SS \subset \TT$ of the big torus.
Let us assume that $X \subset Y$ is a $\SS$-invariant complete intersection 
in $Y$ given by the ideal $I(X) = \langle f_1,\ldots,f_q \rangle$,
where $\deg f_i = \chi_i$.
Then the $\SS$-graded Hilbert series of $X$ is
\[
h_{X, \SS} (t)= h_Y(t) \cdot (1-t^{\chi_1})\cdot \ldots \cdot(1-t^{\chi_q}).
\]
\end{cor}

The next statement is a corollary of Lemma~\ref{hs_of_invariants} 
by using the description of the quotient as the spectrum of invariants given in Theorem~\ref{GIT}. 

\begin{lem}\label{hs_of_git}
Let $X = \Proj R$ be a projective toric variety with 
a subtorus $\iota: \SS \mono \TT$ of the big torus as before.
We assume that both actions are linearized with respect to the ample Weil divisor $\ccL$. 
Let $\iota^* : M_{\TT} \epi M_{\SS}$ be the correspoding surjection of character lattices.
Then
\[
h_{X  \git  \SS} (t_0, t_1,\ldots,t_r) = \sum_{\chi \in \tau(X) \cap \ker \iota^*}
t^{\chi} \cdot \dim R_{\chi}.
\]
\end{lem}
We combine the above facts to get the equality of the Hilbert series 
of models of mutation equivalent graphs.
\begin{thm}\label{same_hs}
Let $\ccG_1$ and $\ccG_2$ be mutation-equivalent graphs and
$S_1$ (respectively $\ccS_2$) be a subset of leaves of $\ccG_1$ 
(respectively $\ccS_2$). Assume that $|\ccS_1|=|\ccS_2|$.
Then the multigraded series with respect to the tori
associated to those sets of leaves are equal
\[
h_{X(\ccG_1),\TT(\ccS_1)} = h_{X(\ccG_2),\TT(\ccS_2)}
\]
\end{thm}
\begin{prf}
Since they are mutation-equivalent by Lemma~\ref{graphs_mutation_equivalent}
we can assume the sequence of mutation takes the set $\ccS_1$
to the set $\ccS_2$. 
We can assume $\ccG_1$ and $\ccG_2$ differ by one mutation.
In Section~\ref{flat_ci_family} 
we have constructed a flat family which is 
a complete intersection having the models $X(\ccG_1)$ and $X(\ccG_2)$ as fibers.
Because both those models are complete intersections of the same type 
in the same ambient space by using Lemma~\ref{hs_of_ci} we conclude 
that the Hilbert series are equal.
\end{prf}
We illustrate Theorem~\ref{same_hs} on examples.
\begin{ex}
We compute the Hilbert series for the models of graphs with 
two leaves and the first Betti number one
$X(\LMfont)$ and $X(\HAMMOCKfont)$.
As we saw in Example~\ref{LMgens} and in the notation of 
Example~\ref{LMembedding}, 
the cone $\tau(\LMfont)$ has generators of the following multidegrees
\[
\begin{array}{|c|c|c|c| c|c|c|c|} \hline
\text{coordinate}
&x_{0000} 
&x_{1100}
&x_{0011}
&x_{1111}
&y_1
&z_1
&z_2\\
\hline
\text{degree}
&(1,0,0)
&(1,0,0)
&(1,1,1)
&(1,1,1)
&(2,1,1)
&(2,2,0)
&(2,0,2)\\\hline
\end{array}
\]
and that $X(\LMfont)$ is a complete intersection in $\PP(1^4,2^3)$ of
a quadric $x_{1100}x_{0011}-x_{0000}x_{1111}$ of multidegree $(2,1,1)$
and a quartic $y_1^2 -z_1z_2$ of multidegree $(4,2,2)$. Thus
\[
h_{X(\LMscriptfont), \ccL(\LMscriptfont)}(t,s_1,s_2)=
\frac
{\cancel{(1-t^2 s_1 s_2)} \cdot (1-t^4 s_1^2 s_2^2)}
{(1-t)^2 \cdot (1-t s_1 s_2)^2 \cdot \cancel{(1-t^2 s_1 s_2)} \cdot (1-t^2 s_1^2) \cdot (1-t^2 s_2^2)}
\]
Ignoring the multigrading by the two dimensional torus spanned by leaves
(by setting $s_1=s_2=1$) we get:
\[
h_{X(\LMscriptfont)}(t)=
\frac
{(1-t^2) \cdot (1-t^4)}
{(1-t)^4 \cdot (1-t^2)^3}
=
\frac
{1+t^2}
{(1-t)^4 \cdot (1-t^2)}
\]
We have seen in Example~\ref{HAMMOCKembedding} that
$\tau(\HAMMOCKfont)$ has six generators: 
those are the elements in the above table for 
$\LMfont$ apart from $y_1$.
The model $X(\HAMMOCKfont)$ is a hypersurface of degree $(4,2,2)$ 
in~$\PP(1^4,2^2)$ so
\begin{multline*}
h_{X(\HAMMOCKfont), \ccL(\HAMMOCKfont)}(t,s_1,s_2)=
\frac
{(1-t^4 s_1^2 s_2^2)}
{(1-t)^2 \cdot (1-t s_1 s_2)^2 \cdot (1-t^2 s_1^2) \cdot (1-t^2 s_2^2)}\\
=h_{X(\LMscriptfont), \ccL(\LMscriptfont)}(t,s_1,s_2).
\end{multline*}

Again  we can ignore the multigrading and get
\[
h_{X(\HAMMOCKfont)}(t)=
\frac
{(1-t^4)} 
{(1-t)^4 \cdot (1-t^2)^2}
=
\frac
{1+t^2}
{(1-t)^4 \cdot (1-t^2)}.
\]

We expand to see the first few terms
\[
h_{X(\LMscriptfont)}(t)=
1 + 4 t + 12 t^2  + 28 t^3  + 57 t^4  + 104 t^5  + 176 t^6  + 280 t^7  + O(t^8).
\]
\end{ex}

\begin{ex}
The Hilbert series of models of both graphs $\THETAfont$ and $\DUMBBELLfont$
with no leaves (thus no additional grading) and two cycles is
\begin{multline*}
h_{X(\DUMBBELLscriptfont )}(t) =
h_{X(\THETAscriptfont) }(t)=
\frac
1
{(t^4 - 4t^3 + 6t^2 - 4t + 1)} = 
1 + 4 t + 10 t^2 + 20 t^3  + 35 t^4  + 56 t^5  \\+ 84 t^6  + 120 t^7  +  O(t^8).
\end{multline*}

This is because $X(\THETAfont)$ is
$\PP^3=(\PP^3\times \PP^3)  \git  (\CC^* \times \CC^* \times \CC^*$).
\end{ex}
\subsection{Computing the Hilbert function.}\label{section_computing_Hilbert_function}

Given a trivalent tree $\ccT$ with $n$ leaves we computed 
the Hilbert function $H_{X(\ccT)}$ of its model 
in \cite{buczynska_wisniewski} as
\[
H_{X(\ccT), \SS(\{l\})}(m,k)=1_m^{\star n}(k),
\]
where the additional grading corresponds to a distinguished leaf $l$,
$\star$ is an appropriate summing formula and 
$1_m$ is the constant function.
This inductive formula for $H_{X(\ccT)}$ uses the decomposition 
of the tree $\ccT$ as a sum of tripods, which leads to the presentation 
of the polytope $\Delta(\ccT)$ as a fiber product of tetrahedrons 
$\Delta(\TRIPOD)$.

The same method works for any trivalent graph.
We proved in Theorem~\ref{same_hs} that the Hilbert function of 
mutation-equivalent graphs are equal. 
By Lemma~\ref{graphs_mutation_equivalent} we know that any  graph 
is mutation-equivalent to a graph of the shape
depicted on Figure~\ref{pic_caterpillar_graph}.

\myVSPACEfigure\begin{center}
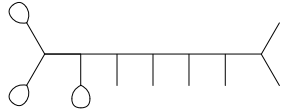
  \begin{minipage}{\textwidth}
\[
\begin{xy}
=(  8, 0     )"A" ;
=(  4, 6.9282)"B" ;
=( -4, 6.9282)"C" ;
=( -8, 0     )"D" ;
=( -4,-6.9282)"E" ;
=(  4,-6.9282)"F" ;
=( 16, 0)"Aleg";
=(-16, 0     )"Dleg";
=(-20, 6.9282)"Hone";
=(-20,-6.9282)"Htwo";
=(-19.6, 7.2282)"HoneXbase";
=(-20.4,-6.6282)"HtwoXbase";
(0,0);(.6,0):;
\save
\POS{"Hone"};\POS{"HoneXbase"}:0*{\CHERRY};
\restore \save
\POS{"Htwo"};\POS{"HtwoXbase"}:0*{\CHERRY};
\restore \save
"Hone"; "Dleg" **\dir{-};
"Htwo"; "Dleg" **\dir{-};
"Dleg" ; "D" **\dir{-};
"Dleg"; (32,0) **\dir{-};
"D" ;   (-8,-6.9282) **\dir{-};
\restore \save
( -8,-6.9282) ; (-8.5, -6.9282):0*{\CHERRY};
\restore \save
(0, 0); ( 0,-6.9282) **\dir{-};
(8, 0); ( 8,-6.9282) **\dir{-};
(16,0); (16,-6.9282) **\dir{-};
(24,0); (24,-6.9282) **\dir{-};
(32,0); (36,-6.9282) **\dir{-};
(32,0); (36, 6.9282) **\dir{-};
\restore \save
(0,15)*{};(0,-15)*{};
\restore
\end{xy}
\] 
\captionof{figure}{Caterpillar graph}
\label{pic_caterpillar_graph}
\end{minipage} \end{center}\myVSPACEfigure
This means we have reduced the calculation to this case of caterpillar graphs.
As we have described in Section~\ref{section_graphs}, any graph is presented 
as union of tripods $\TRIPOD$ with identifications. More precisely, 
any trivalent graph is built from $\TRIPOD$ by the operations of 
grafting two graphs and gluing two leaves.

\begin{rmk}
To produce a caterpillar graph $\ccG$ from 
$\mid$'s (leaves) and $\BALLOONfont$'s (leaves with loop) 
using $\star$ and $\supset$, we need the second operation 
only once per graph and only in the case when $\ccG$ has no leaves.
\end{rmk}
\begin{rmk}
On the level of graph models we have 
\[
X(\ccG_1 \star \ccG_2) = 
\left(X(\ccG_1) \times X(\TRIPOD) \times X(\ccG_2) \right)  \git {(\CC^*)^2}
\]  
and
\[
X(\ccG^{l_1}_{l_2} \glue) = X(\ccG)  \git {\CC^*},
\]
where the actions of the tori were described in 
Section~\ref{section_cone_gens}.
\end{rmk}

We give a formula for 
$h_{X(\ccG_1 \star \ccG_2), \SS(\ccL(\ccG_1\star\ccG_2))}$
and for $h_{X(\ccG^{l_1}_{l_2}\glue), \SS(\ccL(\ccG^{l_1}_{l_2}\glue)}$,
using the above fact about how the model of $\ccG_1 \star \ccG_2$
is constructed from smaller pieces.
\begin{multline}\label{star_formula}
h_{X(\ccG_1 \star \ccG_2), \SS(\ccL(\ccG_1\star\ccG_2))} =
h_{X(\ccG_1), \SS(\ccL(\ccG_1))} \star h_{X(\ccG_2), \SS(\ccL(\ccG_2))} := \\
\text{the part containing monomials of the form } 
(t_1t_2t_3)^i(s')^0(s'')^0 (s''')^j s^I \text{of }\\
h_{X(\ccG_1)}(t_1,s_1,\ldots,s_{n_1},{\frac{1}{s'}}) \cdot 
h_{X(\TRIPOD)}(t_3,s',s'',s''') \cdot
h_{X(\ccG_2)} (t_2,s_{n_1+1},\ldots,s_{n_1+n_2},{\frac{1}{s''}})  
\end{multline}
where $s=(s_1,\ldots,s_{n_2})$ and $I$ is the exponent vector.

Let us compute the input functions: 
apart from the constant one which corresponds 
to leaves of $\ccG$ we have $H_{X(\BALLOONscriptfont), \SS(\{l\})}$ 
the Hilbert function of the model of graph with two edges.
Recall that the model $X(\BALLOONfont)$  is  $\PP^3 \git \CC^* = \PP(1,1,2)$, 
where the $\CC^*$ action has weights $[0 1 0 -1]$ on $\PP^3$. 
Here is the list of generators with weights and 
the resulting graded Hilbert function
\[
\begin{array}{cc}
\begin{array}{|c|c|c|c|}\hline
&\BALLOON..& \BALLOON-.&\BALLOON-=\\\hline
t&1&1&2\\\hline
s&0&0&2\\ \hline
\end{array}&
\leadsto h(t,s)= \frac {1}{(1-t)(1-s^2 t^2)}
\end{array}
\]

We can expand Formula~\eqref{star_formula}, setting $f$ to be a Hilbert function of some graph,
 to get for $k \le \frac {m} {2}$
\begin{align*}
h_{X(\BALLOONscriptfont)} \star f (k) =
(m-k+1) &\sum_{a_0=0}^{m-k-1}  f(a) [2|k+a] (a+1) +\\
(k+1) &\sum_{a_0=m-k}^k  f(a) [2|k+a] (m+2-a) +\\
k   &\sum_{a_0=k+1}^m f(a) [2|k+a] (m+1-a)
\end{align*}
and for $k \ge \frac {m} {2}$
\begin{align*}
h_{\BALLOONscriptfont}\star f (k) =
(m-k+1) &\sum_{a_0=0}^{m-k-1}  f(a) [2|k+a] (a+1) +\\
(m-k+1) &\sum_{a_0=m-k}^k  f(a) [2|k+a] (2m-2k-a+1) +\\
(m-k)   &\sum_{a_0=k+1}^m f(a) [2|k+a] (2m-2k-a)
\end{align*}
where $a=(a_0,a_1, \ldots, a_n)$.

In the same way we can write
\begin{multline}\label{join_formula}
h_{X(\ccG^{l_1}_{l_2} \glue)} (t,s_1,\ldots,s_n) = 
\text{the part that contains monomials } t^i(s')^0s^I \text{ of }
h_{X(\ccG)} (t, s_1, \ldots, s_n, s', {\frac {1}{s'}})
\end{multline}
where $s=(s_1,\ldots,s_n)$ and $I$ is the exponent vector.

\end{document}